\documentclass[aos,MSNbibl,seceqn,nameyear,dvips]{arximspdf}
\usepackage{graphicx}

\doi{10.1214/13-AOS1166} 
\volume{41}
\issue{5}
\pubyear{2013}
\firstpage{2668}
\lastpage{2697}

\makeatletter
\newcommand{\rrVert}{\Vert}
\newcommand{\llVert}{\Vert}
\newtheorem{theorem}{Theorem}
\newtheorem{lemma}{Lemma}
\newtheorem{corollary}{Corollary}
\newtheorem{proposition}{Proposition}
\newproclaim{remark}{Remark}
\newproclaim{example}{Example}
\newproclaim{definition}{Definition}
\makeatother

\begin{document}
\begin{frontmatter}

\title{Spatially inhomogeneous linear inverse problems with possible singularities}
\runtitle{spatially inhomogeneous linear inverse problems}

\begin{aug}
\author[A]{\fnms{Marianna} \snm{Pensky}\corref{}\ead[label=e1]{Marianna.Pensky@ucf.edu}\thanksref{t1}}
\runauthor{M. Pensky}
\affiliation{University of Central Florida}
\address[A]{Department of Mathematics\\
University of Central Florida\\
Orlando, Florida 32816-1354\\
USA\\
\printead{e1}} 
\end{aug}
\thankstext{t1}{Supported in part by NSF Grant DMS-11-06564.}

\received{\smonth{10} \syear{2012}}
\revised{\smonth{8} \syear{2013}}

%
\begin{abstract}
The objective of the present paper is to introduce the concept of a
spatially inhomogeneous linear inverse problem
where the degree of ill-posedness of operator $Q$ depends not only on
the scale but also on location. In this case, the rates of convergence
are determined by the interaction of four parameters, the smoothness
and spatial homogeneity of the unknown function~$f$ and degrees of
ill-posedness and spatial inhomogeneity of operator~$Q$.

Estimators obtained in the paper are based either on wavelet--vaguelette
decomposition (if the norms of all vaguelettes are finite) or on a
hybrid of wavelet--vaguelette decomposition and Galerkin method (if
vaguelettes in the neighborhood of the singularity point have infinite
norms). The hybrid estimator is a combination of a linear part in the
vicinity of the singularity point and the nonlinear block thresholding
wavelet estimator elsewhere. To attain adaptivity, an optimal
resolution level for the linear, singularity affected, portion of the
estimator is obtained using Lepski [\textit{Theory Probab. Appl.}
\textbf{35} (1990) 454--466 and \textbf{36} (1991) 682--697]
method and is used subsequently as the lowest resolution level
for the nonlinear wavelet estimator.
We show that convergence rates of the hybrid estimator lie within a
logarithmic factor of the optimal minimax convergence rates.

The theory presented in the paper is supplemented by examples of
deconvolution with a spatially inhomogeneous kernel and deconvolution
in the presence of locally extreme noise or extremely inhomogeneous
design. The first two problems are examined via a limited simulation
study which demonstrates advantages of the hybrid estimator when the
degree of spatial inhomogeneity is high. In addition, we apply the
technique to recovery of a \mbox{convolution} signal transmitted via amplitude
modulation.
\end{abstract}

%

\end{frontmatter}

\section{Introduction}\label{introduction}
\subsection{Formulation}

Let $Q$ be a known linear operator on a Hilbert space $H$ with inner
product $\langle\cdot, \cdot\rangle$. The objective is to recover $f
\in H$ by observing
%
%
\begin{equation}
\label{eqmaineq} y(x) = (Qf) (x) + \sqrt{\varepsilon} W(x),\qquad x \in
\mathcal{X},
\end{equation}
where $W(x)$ is the white noise process and $\sqrt{\varepsilon}$ is
noise level. Assume that observations can be taken as functionals of
$y$
%
%
\begin{equation}
\label{eqobserv} \langle y, g \rangle= \langle Qf, g \rangle+ \sqrt {\varepsilon
} \xi(g),\qquad g \in H,
\end{equation}
where $\xi(g)$ is a Gaussian random variable with zero mean and
variance $\|g\|^2$ such that $\mathbb{E}[\xi(g_1) \xi(g_2)] = \langle
g_1, g_2 \rangle$. In what follows, $\| \cdot\|$ denotes the
$L^2$-norm, all other norms are explicitly marked.

Model (\ref{eqmaineq}) is a common representation of a linear inverse
problem with the Gaussian noise and has been studied by many authors
[see, \citet{abram}, \citet{bissantz},
\citet{cavalgol2}, \citet{cavalgol1}, \citet{cohen},
\citeauthor{hoffmann} (\citeyear{hoffmann}), \citet{donoho}, \citet{gol},
\citet{kalifa} and \citet{mair},
among others]. A typical assumption in the problem above is that
operator $Q$ acts uniformly over the spaces of functions represented at
a common scale, independently of the location of a function. In
particular, consider a~set of ``test'' functions $ \psi_{ha}(x) =
h^{-1/2} \psi((x-a)/h )$ where $\psi(x), x \in[0,1]$, has a~bounded
support $({L_\psi}, {U_\psi })$ and unit $L^2$-norm $\|\psi\| =1$.
Then, functions $\psi_{ha}(x)$ have scale $h$, supports concentrated
around $x=a$ and unit norms. Conditions which are commonly imposed on
operator $Q$ imply that it contracts the norms of all functions
$\psi_{ha}$ uniformly, that is, the value of $\| Q \psi_{ha}\|$ depends
considerably on the scale $h$ but hardly at all on $a$. Moreover, if
there exist $(Q^*)^{-1} \psi_{ha}$, where $Q^*$ is the adjoint of
operator $Q$, then values of $\|(Q^*)^{-1} \psi_{ha}\|$ follow the same
pattern. However, not all linear operators necessarily have those
properties.

In order to illustrate the discussion above, consider linear operator
$Q$ with the adjoint $Q^*$ given by
%
%
\begin{equation}
\label{eqQQstar} (Qf) (x) = \mu(x) \int_0^x
f(t) \,dt,\qquad \bigl(Q^* v \bigr) (x) = \int_x^1
\mu(z) v(z) \,dz,
\end{equation}
where $\mu(x)$ is a smooth function. Assume that function $\psi$ 
is continuously differentiable and integrates to zero: $\int_0^1
\psi(z) \,dz =0$. Denote $\Psi(z) = \int_{L_\psi}^z \psi(x) \,dx$ and
observe that $\Psi(z) =0$ whenever $z \notin({L_\psi}, {U_\psi})$.
Then, direct calculations yield
$ (Q \psi_{ha}) (y) = h^{1/2} \mu(y) \Psi(\frac{y-a}{h} )$,  $
(Q^*)^{-1} \psi_{ha}(y) = - h^{-3/2} \mu^{-1} (y) \psi^\prime(
\frac{y-a}{h} )$,  so\vspace*{-1pt} that $\| Q \psi_{ha}\|^2 = h^2 \int_{L_\psi}
^{U_\psi} \mu^2(a + hz) \Psi^2 (z) \,dz = h^2 [\mu^2 (a) \| \Psi\|^2 +
o(1) ]$ as \mbox{$h \to0$} and $\|(Q^*)^{-1} \psi_{ha}\|^2 = h^{-2}
\int_{L_\psi}^{U_\psi}\mu^{-2}(a + hz) [\psi^\prime(z)]^2 \,dz. $

If $\mu(y)$ is a constant or, at least, $C_\mu^{-1} < \mu(y) < C_\mu$
for some relatively small~$C_\mu$, then dependence of $\| Q \psi
_{ha}\|$ and $\|(Q^*)^{-1} \psi_{ha}\|$ on $a$ can be ignored, so
equation~(\ref{eqmaineq}) with $Q$ given by (\ref{eqQQstar}) can be
treated as a spatially homogeneous problem. However, if $C_\mu$ is
large, dependence on $a$ becomes essential and \mbox{equation
(\ref{eqmaineq})} is a spatially inhomogeneous inverse problem.

Dependence on $a$ becomes even more extreme if $\mu(y)$ vanishes at
some point $x_0 \in(0,1)$, for example, $\mu^2 (x) = C_\alpha|x -
x_0|^\alpha$. Indeed, in this case, $x_0$ is the \textit{singularity
point} and it is easy to show that $\| Q \psi_{hx_0} \|^2 \asymp
h^{2+\alpha}$ and
\[
\bigl\| \bigl(Q^* \bigr)^{-1} \psi_{hx_0} \bigr\|^2 \asymp
\cases{h^{-(2 + \alpha)}, &\quad if $\alpha<1$, \vspace*{2pt}
\cr
\infty, &\quad if
$\alpha\geq1$.}
\]

Since wavelets provide an adequate tool for scale-location
representations of functional spaces, it is convenient to introduce
spatially inhomogeneous linear inverse problems using a
wavelet--vaguelette decomposition proposed by \citet{donoho}. In
particular, in the case when $H= L^2(\mathcal{D})$, $\mathcal{D}
\subset R$, Donoho's assumptions appear as follows:
\begin{longlist}[(D3)]
\item[(D1)] There exist three sets of functions: $\{
\psi_{jk}\}$,
an orthonormal wavelet basis of $H$, and nearly orthogonal sets $\{
u_{jk} \}$ and $\{ v_{jk} \}$ such that $ Q \psi_{jk} = v_{jk}$, $Q^*
u_{jk} = \psi_{jk}$, $\|v_{jk} \| = \lambda_{j}$, $\|u_{jk} \| \asymp
\lambda_{j}^{-1}$,  where $\lambda_{j}$ depend on resolution index
$j$ but not on spatial index $k$.

\item[(D2)] $u_{jk}$ and $v_{jk}$ are such that $\langle u_{j_1
k_1}$, $v_{j_2 k_2} \rangle= \delta_{j_1, j_2} \delta_{k_1, k_2}$.

\item[(D3)] Sets $ \{u_{jk} \}$ and $\{ v_{jk} \}$ are nearly
orthogonal, that is, for any sequence $\{a_{jk} \} \in l^2$ one has
%
\[
\biggl\llVert\sum_{j,k} a_{jk}
\lambda_{j} u_{jk} \biggr\rrVert^2 \asymp\sum
_{j,k} a_{jk}^2,\qquad \biggl
\llVert \sum_{j,k} a_{jk}
\lambda_{j}^{-1} v_{jk} \biggr
\rrVert^2 \asymp\sum_{j,k}
a_{jk}^2.
\]
Under conditions (D1)--(D3), $f$ can be recovered using reproducing
formula
%
%
\begin{equation}
\label{eqrepr1} f = \sum_{j,k} \langle Qf,
u_{jk} \rangle\psi_{jk},
\end{equation}
which is analogous to the reproducing formula for the SVD. Assumptions
\mbox{(D1)--(D3)} are quite standard. Indeed, similar assumptions were
introduced in \citet{cavalgol1}, \citet{cavalgol2},
\citet{gol} and \citet{knapik}. The common premise is that
operator $Q$ acts ``uniformly'' over subspaces of $H$, so singular
values or their surrogate equivalents depend on the resolution level
only but not on location. If $V_j = \operatorname{Span}\{
\psi_{i,k}, i \leq j, k \in\mathbb{Z}\}$ is the subspace of
functions at resolution level $j$, the above assumptions reduce to a
common assumption of Galerkin method [see, e.g., \citet{cohen} or
\citet{hoffmann}] that on subspace $V_j$\vadjust{\goodbreak} operator $Q$ has a
bounded inverse with the norm dependent on $j$ only, that is, there
exist $\lambda_j>0$ such that
%
%
\begin{equation}
\label{eqGelerkinass} \sup_{j \geq
0} \bigl[ \lambda_j
\bigl\|Q^{-1} \bigr\|_{V_j \rightarrow Q^{-1}(V_j)} \bigr]< \infty,
\end{equation}
which is very similar to combination of assumptions (D1) and (D3)
above.

Note that both, assumptions (D1) and (\ref{eqGelerkinass}) imply that
any function $v \in V_j$ with $\|v\|=1$ has an inverse image, the norm
of which is bounded by a constant which is independent of the support
of $v$. In this sense, operator $Q$ is an ill-posed \textit{spatially
homogeneous} operator. In the present paper, we shall be interested in
a~different situation when assumptions (D1) and (D3) may not be true.
In particular, we assume that the norms of the inverse images of
$\psi_{j,k}$ depend on the spatial index $k$ and may be unbounded, that
is, condition (D1) and possibly condition (D3) are violated. We shall
refer to the such inverse linear problems as \textit{spatially
inhomogeneous} in comparison with spatially homogeneous problems which
satisfy conditions (D1)--(D3) above.
\end{longlist}

\subsection{Motivation}

Spatially inhomogeneous ill-posed problems appear naturally in the case
when either the noise level is spatially dependent or observations are
irregularly spaced. Problems of this kind have been considered
previously, both theoretically and in practical applications.
Nevertheless, in former studies, it was always assumed that the noise
level is uniformly bounded above or the design density of observations
is bounded away from zero. The situations investigated in the present
paper rather refer to \textit{locally extreme noise} or
\textit{extremely inhomogeneous design} (which can be also described as
\textit{a local data loss}). Traditionally, in the first situation,
measurements are treated as outliers and are removed from future
analysis while the second one is dealt with using missing data
techniques. Approach suggested in the present paper provides an
alternative to those methodologies.
Extreme noise or extremely inhomogeneous design occur in analysis of
forensic data [see, e.g., \citet{li}] and deconvolution of LIDAR
signals [see, e.g., \citet{harsdorf} and \citet{gurdev}]
or astronomical images [see, e.g., \citet{starck}]. 
%
In addition, spatially inhomogeneous ill-posed problems arise whenever
the kernel is spatially inhomogeneous, as, for example, in the case of the
amplitude modulation.
%
Below we consider some examples in more detail.

%
%
\begin{example}[(Deconvolution of LIDAR signals)]\label{ex1}
LIDAR (Light Detection And Ranging or Laser Imaging Detection And
Ranging) is an optical remote sensing technology that can measure the
distance to, or other properties of, targets by illuminating the target
with laser light and analyzing the backscattered light. LIDAR
technology has applications in archaeology, geography, geology,
geomorphology, seismology, forestry, remote sensing, atmospheric
physics. LIDAR data model is mathematically described by convolution\vadjust{\goodbreak}
equation $P=R*P_\delta$ where $P$ is the time-resolved LIDAR signal,
$P_\delta$ is the impulse response function and $R$~is the system
response function to be determined [see, e.g., \citet{harsdorf}
and \citet{gurdev}]. However, if the system response function of
the LIDAR is longer than the time resolution interval, then the
measured LIDAR signal is blurred and the effective accuracy of the
LIDAR decreases. This loss of precision becomes extreme when, for
example, LIDAR is used to for emergency response and natural disaster
management such as assessment of the extent of damage due to volcanic
eruptions or forest fires.
In this situation, routinely, distances are calculated through
filtering of the data set (removing outliers) and applying
interpolation techniques. However, keeping all existing data and
accounting for extreme noise may improve precision of the analysis of
LIDAR signals.
\end{example}

\begin{example}[(Amplitude modulation)] \label{ex3}
Amplitude Modulation (AM) is a way of transmitting information in the
form of electro-magnetic waves. In AM, a radio wave known as the
``carrier'' is modulated in amplitude by the signal, that is, to be
transmitted, while the frequency remains constant [see, e.g.,
\citet{miller}]. In video or image transmission (such as TV) where
the base-band signal has inherent large bandwidth, AM is usually
preferred to Frequency Modulation (FM) systems since the latter ones
require additional bandwidth. Since in an AM, signal information is
``stored'' in amplitude which is affected by noise, AM is more
susceptible to noise than FM. Mathematically, the problem reduces to
multiplying the transmitted signal by the function $\mu(x) =
\cos(2\pi\omega x - \theta)$ with large $\omega\approx n/2$ and
$\theta\in[0; 2\pi]$. In Section~\ref{secrealdata}, we provide an
in-depth description of application of the methodology developed in the
paper to recovery of a~convolution signal transmitted via AM.
\end{example}

\subsection{Objectives and layout of the paper}

The objective of the present paper is to introduce the concept of a
spatially inhomogeneous linear inverse problem which, to the best of
the author's knowledge, has never been considered previously in
statistical framework. It turns out that spatially inhomogeneous
problems exhibit properties which are very different from their
spatially homogeneous counterparts. In particular, if the norms of
vaguelettes $u_{j k}= (Q^*)^{-1} \psi_{jk}$ are infinite in the
vicinity of a singularity point, reproducing formula (\ref{eqrepr1})
cease working and the usual wavelet--vaguelette estimators cannot be
applied. In this case, we propose a hybrid estimator which is based on
combination of wavelet--vaguelette decomposition and Galerkin method. We
study two application of the general theory, deconvolution with
spatially inhomogeneous design and deconvolution with a spatially
inhomogeneous kernel (the case of heterogeneous noise being a
particular case of the latter).

Another interesting feature of the model is that the rates of
convergence are determined by the interaction of four parameters, the
smoothness and spatial homogeneity of the unknown function $f$ and
degrees of ill-posedness and spatial inhomogeneity of operator $Q$. In
particular, if operator $Q$ is weakly inhomogeneous, then the rates of
convergence are not influenced by spatial inhomogeneity of operator $Q$
and coincide with the rates which are usual for homogeneous linear
inverse problems.

In what follows, we assume that operator $Q$ in (\ref{eqmaineq}) is
completely known. If, in practical applications, this is not true, one
has to account for the extraneous errors which stem from the
uncertainty in the operator $Q$ by using, for example, ideas of
\citet{hoffmann}. Also, to simplify our considerations, we limit our
study to the case when $\mathcal{X} = [0,1]$, $H= L^2[0,1]$ and $k$ is
a scalar. The theory presented below can be generalized to the case
when $H= L^2(\mathcal{D})$, $\mathcal{D} \subset R^d$ and $k$ is a
$d$-dimensional vector. This extension should be relatively
straightforward if one is dealing with isotropic Besov spaces but
becomes much more interesting and involved in the case of anisotropic
Besov spaces [see, e.g., \citet{kerk}]. However, we leave those
extensions for future investigations since considering them below will
prevent us from focusing on the main objective of the paper.


The rest of the paper is organized as follows.
Section~\ref{secassumptions} introduces the concept of a spatially
inhomogeneous ill-posed problem and formulates major definitions and
assumptions which are used throughout the paper.
Section~\ref{seclowbounds} presents the asymptotic minimax lower bounds
for the $L^2$-risk of the estimators of the solution of the problem
over a wide range of Besov balls. Section~\ref{secestimationstrategies}
talks about estimation strategies, in particular, about partitioning
the unknown response function $f$ and its estimator into the
singularity-affected and the singularity-free parts, the main idea at
the core of the hybrid estimator. Section~\ref{secriskfcf0} elaborates
on the risk of the estimator constructed in the previous section when
the lowest resolution level in the zero-affected portion of the
estimator is fixed. Section~\ref{secadaptive} discusses the adaptive
choice of the lowest resolution level resolution level and derives the
asymptotic minimax upper bounds for the $L^2$-risk. In
Section~\ref{secexamples}, we consider two examples of spatially
inhomogeneous ill-posed problems, deconvolution with the spatially
inhomogeneous operator (Section~\ref{secdeconv}) which can be viewed as
a version of a deconvolution equation with spatially inhomogeneous
noise, and deconvolution based on irregularly spaced sample
(Section~\ref{secirreg}). Section~\ref{secsimulations} presents a
limited simulation study of deconvolution with heteroscedastic noise
and also studies application of the hybrid estimator to recovery of a
convolution signal transmitted via amplitude modulation.
Section~\ref{secdiscussion} concludes the paper with a discussion.
Proofs of the statements are contained in the supplementary material [\citet{pensuppl}].


\section{Spatially inhomogeneous ill-posed problem: Assumptions and
definitions} \label{secassumptions}

Consider a scaling function $\varphi$ and a corresponding wavelet
$\psi$ with bounded supports and form an orthonormal\vadjust{\goodbreak} wavelet basis $\{
\psi_{jk}\}$ of $L^2([0,1])$. We further impose the following set of
assumptions on spatially inhomogeneous operator $Q$.
\begin{longlist}[(A3)]
\item[(A1)] There exist functions $\{ u_{j k} \}$ and $\{ v_{j k}\}$
such that $Q \psi_{jk}= v_{j k}$, $Q^* u_{j k}= \psi_{jk}$,
where $\| v_{j k}\| = \lambda_{j k}< \infty$.

\item[(A2)] There exists a \textit{singularity} point $x_0\in(0,1)$
and a constant $D \geq0$ such that $\|u_{j k}\| = \infty$ if $ |k -
k_{0j}| < D$ and, for any $j$ and any $\{a_{j k}\}$, $k=0,\ldots,
2^j-1$, one has
%
%
\begin{equation}
\label{eqill-posedness} \biggl\llVert\sum_{|k - k_{0j}|
\geq D}
a_{j k}\lambda_{j k}u_{j k} \biggr\rrVert
^2 \leq C_u \sum_{|k - k_{0j}|
\geq D}
a_{j k}^2,
\end{equation}
where $C_u<\infty$ is independent of $j$ and $k_{0j}= 2^j x_0$ is the
parameter corresponding to location $x_0$ ($k_{0j}$ is not necessarily
an integer).

\item[(A3)] Functions $v_{j k}$ are such that, for any $j$,
inequality
%
%
\begin{equation}
\label{equpperunifv} \Biggl\llVert\sum_{k=0}^{2^j -1}
a_{j k} \lambda_{j k}^{-1} v_{j k} \Biggr
\rrVert^2 \leq C_v \sum_{k=0}^{2^j -1}
a_{j k}^2
\end{equation}
holds for any $\{a_{j k}\}$, $k=0,\ldots, 2^j-1$, where $C_v < \infty$
is independent of $j$.

Note that assumptions (A1)--(A3) are weaker than assumptions (D1)--(D3)
above. First, $\lambda_{j k}$ depends not only on resolution level but
also on location of the wavelet coefficient.
Also, if $D>0$, then, in the neighborhood of the singularity point
$x_0$, wavelet coefficients cannot be recovered directly since $\|u_{j
k}\| = \infty$, and we say that operator $Q$ has a \textit{singularity}
at $x_0$.

Since one usually start wavelet expansion at some finite resolution
level $m$, we need an extra assumption which mirrors assumption (A2)
and can be derived from~it:

\item[(A4)] There exist functions $\{ t_{m k} \}$ and positive
constants $C_t$ and $D_0$ independent of $m$ such that, for any
$\{a_{mk}\}$, $k=0,\ldots, 2^m-1$,
%
%
\begin{equation}
\label{eqtjk} \quad Q^* t_{m k}= \varphi_{mk}, \qquad \biggl
\llVert\sum_{|k - k_{0j}| \geq D_0} a_{mk}
\lambda_{m k}t_{m
k} \biggr\rrVert^2 \leq
C_t \sum_{|k - k_{0j}| \geq D_0} a_{mk}^2.
\end{equation}

If $D=D_0 =0$ in assumptions (A2) and (A4), then $\|u_{j k}\| < \infty$
and $\|t_{m k}\| < \infty$ for any $k$. Hence, $f$ can be expressed
using reproducing formula (\ref{eqrepr1}) which, in this case, becomes
%
%
\begin{equation}
\label{funf} f(x) = \sum_{k=0}^{2^{m}-1}
a_{mk}\varphi_{mk} (x) + \sum_{j=m}^\infty
\sum_{k=0}^{2^{j}-1} b_{jk}
\psi_{jk}(x),
\end{equation}
where $a_{mk}= \langle Qf, t_{m k}\rangle$ and $b_{jk}= \langle Qf,
u_{j k}\rangle$. If $D>0$, reproducing formula (\ref{funf}) cease
working and one needs an alternative solution to recovering $f$.
Indeed, if $Qf$ in expressions for $a_{mk}$ and $b_{jk}$ is replaced by
$y = Qf + \varepsilon W$, then the variances of the wavelet
coefficients in the vicinity of singularity $x_0$ are infinite:
\mbox{$\operatorname{Var}\langle y, u_{j k}\rangle= \infty$} if $ |k -
k_{0j}| < D$ and similar consideration applies to $t_{m k}$. For this
reason, at each resolution level, we partition the set of all indices
into the \textit{singularity-affected} indices
\begin{eqnarray*}
K_{0m} & = & \bigl\{ k=0,\ldots, 2^m-1\dvtx|k -
k_{0m}| < D_0 \bigr\},
\\
K_{1j} & = & \bigl\{ k=0,\ldots, 2^j-1\dvtx|k -
k_{0j}| < D \bigr\}
\end{eqnarray*}
and the \textit{singularity-free} indices
\begin{eqnarray*}
K_{0m}^c&=& \bigl\{k\dvtx0 \leq k \leq2^m -1,
k \notin K_{0m} \bigr\},
\\
K_{1j}^c&=& \bigl\{k\dvtx0 \leq k \leq2^j -1,
k \notin K_{1j} \bigr\}.
\end{eqnarray*}

To be specific, in what follows, we assume that $\lambda_{j k}$ are
such that, for some positive constants $\alpha, \beta$, $C_{\lambda_0}$
and $C_\lambda$ independent of $j$ and $k$, one has
%
%
\begin{equation}
\label{eqlambdajk} \quad C_{\lambda_0} 2^{-j(\alpha+ \beta)} \bigl(1 + |k -
k_{0j}|\bigr)^\alpha\leq\lambda_{j k}^2 \leq
C_\lambda2^{-j(\alpha+ \beta)} \bigl(1 + |k - k_{0j}|\bigr)^\alpha.
\end{equation}
We shall refer to coefficients $\beta$ and $\alpha$ in
(\ref{eqlambdajk}) as degrees of \textit{ill-posedness} and
\textit{spatial inhomogeneity}, respectively. Observe that with
$\lambda_{j k}$ satisfying condition (\ref{eqlambdajk}), the variances
of the coefficients at the lower resolution levels may be significantly
higher than the variances of the coefficients at higher resolution
levels as long as the locations of the lower resolution level
coefficients lie in a close proximity of a singularity point.

In the present paper, we consider estimation of a solution of
inhomogeneous linear inverse problems in the case when the unknown
function $f$ is possibly spatially inhomogeneous itself, in particular,
$f$ belongs to a Besov ball $B_{p,q}^s (A)$ of radius $A$. Interplay
between spatial inhomogeneity of operator $Q$ and properties of~$f$
lead to various very interesting phenomena. In particular, if $\alpha$
is small or $p$ and $\beta$ are relatively large, spatial inhomogeneity
does not affect convergence rates and $f$ can be recovered as well as
in the case of $\alpha=0$.
\end{longlist}

%
%
\begin{remark}[(Multiple singularity points)]\label{remseveralsing}
Note that one can consider a spatially inhomogeneous problems with
multiple singularity points $x_{0,1}< x_{0,2} < \cdots< x_{0,L}$ and
corresponding constants $D_1,\ldots, D_L$ where $L<\infty$ and $x_{0,i}
- x_{0,i-1} \geq\delta >0$ for some fixed positive $\delta$. The theory
developed below can be easily extended to this case, with the
convergence rates of the estimators determined by the ``worst case
scenario'' among singular points $x_{0,i}$, $i=1,\ldots, L$.
\end{remark}

\section{Minimax lower bounds for the risk over Besov balls}
\label{seclowbounds}

Before constructing an estimator of the unknown function $f$ under
model (\ref{eqmaineq}), we derive the asymptotic minimax lower bounds
for the $L^2$-risk over a wide range of Besov balls.

Recall that for an $r_0$-regular multiresolution analysis [see, e.g.,
\citet{meyer1992}, pp. 21--25], with $0<s<r_0$, and for a Besov
ball $ B_{p,q}^s (A) =\{f \in L^p([0,1])\dvtx f \in B_{p,q}^s,
\|f\|_{B_{p,q}^s} \leq A \}, $ of radius $A>0$ with $1\leq p$, $q \leq
\infty$ and $s' = s+1/2-1/p$, one has
%
\[
\|f\|_{B_{p,q}^s} = \Biggl(\sum_{k=0}^{2^m -1}
|a_{mk}|^p \Biggr)^{1/p} + \Biggl(\sum
_{j=m}^{\infty} 2^{js'q} \Biggl(\sum
_{k=0}^{2^j -1} |b_{jk}|^p
\Biggr)^{q/p} \Biggr)^{1/q}
\]
%
%
with respective sum(s) replaced by maximum if $p=\infty$ and/or
$q=\infty$ [see, e.g., \citet{jpkr}].
%

In what follows, we use the symbol $C$ for a generic positive constant,
which takes different values at different places and is independent of
the noise level $\varepsilon$.
The following statement provides the asymptotic minimax lower bounds
for the $L^2$-risk over Besov balls $B_{p,q}^s (A)$.
%

%
%
\begin{theorem} \label{thlower}
Let $1 \leq p, q \leq\infty$ and $s > \max(1/p, 1/2)$. Then, under
assumptions \textup{(A1)--(A3)}, as $\varepsilon\rightarrow0$,
%
%
\begin{equation}
\label{low} R_\varepsilon \bigl(B_{p,q}^s (A) \bigr) =
\inf_{\tilde{f}} \sup_{f\in B_{p,q}^s (A)} \mathbb{E} \|\tilde{f}
- f\|^2 \geq C \Delta(\varepsilon),
\end{equation}
where the infimum is taken over all possible square-integrable
estimators $\tilde{f} $ of $f$ based on $y$ from model (\ref{eqmaineq})
and
%
%
\begin{equation}\label{eqDDeps}
\Delta(\varepsilon) = \cases{
A^{(2(\alpha+ \beta))/(2s' + \alpha+ \beta)}\varepsilon^{2s'/(2s' + \alpha+ \beta)},
\vspace*{2pt}\cr
\qquad \mbox{if }2 s (\alpha-1) \geq(\beta+1) (1-2/p),
\vspace*{4pt}\cr
A^{(2(\beta+ 1))/(2s + \beta+ 1)} \varepsilon^{2s/(2s +\beta+ 1)},
\vspace*{2pt}\cr
\qquad \mbox{if }2 s (\alpha-1) < (\beta+1) (1-2/p).}
\end{equation}
\end{theorem}

%
%
\begin{remark}[(Convergence rates)]  \label{remratediscuss}
As we show below, the minimax global convergence rates in Theorem
\ref{thlower} are attainable up to a logarithmic factor. The rates are
determined by the interaction of four parameters, $s, p, \alpha$ and
$\beta$. Parameters $s$ and $p$ describe, respectively, smoothness and
spatial homogeneity of the unknown function $f$, while $\beta$ and
$\alpha$, defined in (\ref{eqlambdajk}), are referred to as degrees of
ill-posedness and spatial inhomogeneity of operator $Q$. If the value
of $\alpha$ is large, in particular, $2 s \alpha> 2s' + \beta(1 -
2/p)$, convergence rate is significantly affected by the degree of
spatial inhomogeneity $\alpha$ of $Q$. On the other hand, if $2 s
\alpha< 2s' + \beta(1 - 2/p)$, spatial inhomogeneity of operator $Q$
does not affect convergence rate which is determined entirely by the
degree of ill-posedness $\beta$.
\end{remark}


\section{Estimation strategies in the presence of a singularity}
\label{secestimationstrategies}

To be more specific, consider a periodized version of the wavelet basis
on the unit interval
%
%
\begin{equation}
\label{eqorthbasis} \bigl\{\varphi_{mk}, \psi_{jk}\dvtx j
\geq m, k=0,1,\ldots, 2^{j}-1 \bigr\},\vadjust{\goodbreak}
\end{equation}
where $\varphi_{mk}(x) = 2^{m/2} \varphi(2^m x -k), \psi_{jk}(x) =
2^{j/2} \psi(2^j x -k), x \in[0,1]. $ Note that the latter requires
that the resolution level $m$ is high enough, in particular, $m \geq
m_1$, where $m_1$ is such that
%
%
\begin{equation}
\label{eqm1val} 2^{m_1} > \max(\mathcal{L}_{\varphi^*},
\mathcal{L}_{\psi^*} ).
\end{equation}
Here, $\mathcal{L}_{\varphi^*}$ and $\mathcal{L}_{\psi^*}$ are the
lengths of supports of the mother and father wavelets, $\varphi^*$ and
$\psi^*$, that generate periodized wavelet basis. Then, for any $m \geq
m_1$, the set (\ref{eqorthbasis}) forms an orthonormal wavelet basis
for $L^2([0,1])$ and, hence, any $f\in L^2([0,1])$, can be expanded
using formula (\ref{funf}). Under assumptions (A1), (A2) and (A4), one
can construct unbiased estimators of coefficients $a_{mk}$ and $b_{jk}$
%
%
\begin{equation}
\label{coefest} \hat{a}_{mk}= \langle y, t_{m k}\rangle,
\qquad\hat{b}_{jk}= \langle y, u_{j k}\rangle.
\end{equation}
If $k \in K_{0m}^c$ and $k \in K_{1j}^c$, respectively, then estimators
$\hat{a}_{mk}$ and $\hat{b}_{jk}$ have finite variances
%
%
\begin{eqnarray}\label{eqvariances}
\operatorname{Var} (\hat{a}_{mk}) &\asymp&
\lambda_{m k}^{-2},\qquad k \in K_{0m}^c,
\nonumber\\[-8pt]\\[-8pt]
\operatorname{Var}(\hat{b}_{jk}) &\asymp& \lambda_{j k}^{-2},\qquad
k \in K_{1j}^c\nonumber
\end{eqnarray}
and have infinite variances otherwise. In order to account for the
latter, for any $m \geq m_1$, we partition $f$ into the sum of
singularity-affected and singularity-free parts
\[
f(x) = f_{0,m}(x) + f_{c,m}(x),\qquad x \in[0,1],
\]
where
%
%
\begin{eqnarray}
\label{fo}
\quad f_{0,m}(x) & = & \sum_{k \in K_{0m}}
a_{mk} \varphi_{mk} (x) + \sum
_{j=m}^{\infty} \sum_{k \in K_{1j}}
b_{jk} \psi_{jk}(x),\qquad x \in[0,1],
\\
\label{fc} f_{c,m}(x) & = & \sum_{k \in K_{0m}^c}
a_{mk} \varphi_{mk}(x) + \sum_{j=m}^{\infty}
\sum_{k \in K_{1j}^c} b_{jk} \psi_{jk}(x),
\qquad x \in[0,1].
\end{eqnarray}
We then construct estimators $\hat{f}_{0,m}$ and $\hat{f}_{c,m}$ of
$f_{0,m}$ and $f_{c,m}$, respectively, and estimate $f$ by a hybrid
estimator
%
%
\begin{equation}
\label{totalest} \hat{f}_{m}(x) = \hat{f}_{0,m}(x) +
\hat{f}_{c,m} (x),\qquad x \in[0,1].
\end{equation}
In particular, we shall use a linear estimator with the resolution
level $m$ estimated from the data as $\hat{f}_{0,m}$ and a nonlinear
block thresholding wavelet estimator as $\hat{f}_{c,m}$, where the
lowest resolution level $m$ in $\hat{f}_{c,m}$ is determined by the
linear part $\hat{f}_{0,m}$. In what follows, we shall consider
estimation of $f_{0,m}$ and $f_{c,m}$ separately.


First, we construct a block thresholding wavelet estimator $\hat
{f}_{c,m}$ of $f_{c,m}$. For this purpose, we divide the wavelet
coefficients at each resolution level into $l_{j\varepsilon}^L$ blocks
of length $\ln(\varepsilon^{-1})$ to the left\vadjust{\goodbreak} of $(k_{0j}-D)$ and
$l_{j\varepsilon}^R$ blocks to the right of $(k_{0j}+ D)$, where
%
%
\begin{eqnarray}\label{eqljedef}
l_{j\varepsilon}^L &=& (k_{0j}- D)/\ln \bigl(\varepsilon^{-1} \bigr),\nonumber
\\
l_{j\varepsilon}^R &=& \bigl(2^j - D -1 - k_{0j} \bigr)/\ln \bigl(\varepsilon^{-1} \bigr),
\\
l_{j \varepsilon}&=& \max\bigl(l_{j\varepsilon}^L,l_{j\varepsilon}^R \bigr).\nonumber
\end{eqnarray}
Define blocks $U_{jl}^L$ and $U_{jl}^R$ of indices $k$ to the left of
$(k_{0j}-D)$ and to the right of $(k_{0j}+ D)$, respectively, as
\begin{eqnarray*}
U_{jl}^L&=& \bigl\{k\dvtx-D -l \ln \bigl(
\varepsilon^{-1} \bigr)< k - k_{0j}\leq-D - (l-1) \ln \bigl(
\varepsilon^{-1} \bigr) \bigr\},\qquad l \in U_j^L,
\\
U_{jl}^R&=& \bigl\{k\dvtx D + (l-1) \ln \bigl(
\varepsilon^{-1} \bigr)< k - k_{0j}\leq D + l \ln \bigl(
\varepsilon^{-1} \bigr) \bigr\},\qquad l \in U_j^R,
\end{eqnarray*}
where
%
%
\begin{eqnarray}\label{equjlr}
U_j^L &=& \bigl\{l\dvtx1 \leq l \leq
l_{j\varepsilon}^L \bigr\},\nonumber
\\
U_j^R &=& \bigl\{l\dvtx1 \leq l \leq l_{j\varepsilon}^R \bigr\},
\\
U_{j}&=& U_j^L \cup U_j^R.\nonumber
\end{eqnarray}
To simplify the narrative, we shall write $l \in U_{j}$ and $k \in
U_{jl}$ without a specific reference whether a block lies to the right
or to the left of $k_{0j}$. Denote
%
%
\begin{eqnarray}
\label{eqbjl} B_{jl}& = & \sum_{k \in U_{jl}}b_{jk}^2,
\qquad \widehat{B}_{jl}= \sum_{k \in U_{jl}}
\hat{b}_{jk}^2,
\\
\label{eqrlje} R_{j l \varepsilon}& = & \frac{\varepsilon\ln(\varepsilon^{-1})
2^{j(\alpha+ \beta)}}{ |D + (l-1)
\ln(\varepsilon^{-1})|^\alpha} \asymp\varepsilon\sum
_{k \in U_{jl}}\lambda_{j k}^{-2}.
\end{eqnarray}
For any $m \geq m_1$, estimate $f_{c,m}$ by
%
%
\begin{eqnarray}
\label{hfcmfixed}
\hat{f}_{c,m}(x) &=& \sum
_{k \in K_{0m}^c} \hat{a}_{mk}\varphi_{mk}(x)
\nonumber\\[-8pt]\\[-8pt]
&&{} + \sum_{j=m}^{J-1} \sum
_{l \in U_j} \sum_{k \in U_{jl}}\hat
{b}_{jk}\mathbb{I} \bigl( \widehat{B}_{jl}\geq
\tau^2 R_{j l \varepsilon} \bigr) \psi_{jk}(x),\nonumber
\end{eqnarray}
where $\mathbb{I}(\Omega)$ is the indicator function of the set
$\Omega$, the value of $\tau$ will be defined later and
%
%
\begin{equation}
\label{Jvalueo} 2^J = \varepsilon^{-2/(\alpha+\beta+ 2)}.
\end{equation}


Now, consider estimation of the singularity-affected part. Since the
estimators $\hat{a}_{mk}$ of $a_{mk}$, given in (\ref{coefest}), have
infinite variances when $k \in K_{0m}$, we estimate those coefficients
by solving a system of linear equations. Denote $w_{m k}= Q
\varphi_{mk}$ and observe that, for a given $m$, $m_1 \leq m \leq J-1$,
one has $f = f_m + R_m$. Here
%
%
\begin{eqnarray}
\label{fmRm}
f_m &=& \sum_{k \in K_{0m}}
a_{mk} \varphi_{mk}+ \sum_{k \in K_{0m}^c}
a_{mk} \varphi_{mk},
\nonumber\\[-8pt]\\[-8pt]
R_m &=& \sum
_{j=m}^{\infty} \sum_{k=0}^{2^j-1}
b_{jk} \psi_{jk}\nonumber
\end{eqnarray}
and, hence,
%
%
\begin{equation}
\label{eqsum} Qf = \sum_{k \in K_{0m}} a_{mk}
w_{m k}+ \sum_{k \in K_{0m}^c} a_{mk}
w_{m k}+ Q R_m.
\end{equation}
Taking scalar products of both sides of (\ref{eqsum}) with $w_{m l}$,
$l \in K_{0m}$, obtain
%
%
\begin{eqnarray}
\label{eqscprod}
\qquad\langle w_{m l}, Qf \rangle& =& \sum
_{k \in K_{0m}} a_{mk} \langle w_{m l},
w_{m k} \rangle
\nonumber
\\[-8pt]
\\[-8pt]
&&{} + \sum_{k \in K_{0m}^c} a_{mk} \langle
w_{m l}, w_{m k}\rangle+ \langle w_{m l}, Q
R_m \rangle,\qquad l \in K_{0m}.
\nonumber
\end{eqnarray}
%
%
Introduce matrices $\mathbf{A}^{(m)}$ and $\mathbf{B}^{(m)}$ and
vectors $\mathbf{c}^{(m)}$, $\hat{\mathbf{c}}{}^{(m)}$,
$\mathbf{r}^{(m)}$, $\mathbf{z}^{(m)}$, $\mathbf{h}^{(m)}$ and
$\hat{\mathbf{h}}{}^{(m)}$ with elements
%
%
\begin{eqnarray}
\label{matrAB} A^{(m)}_{lk} & =& \langle w_{m l},
w_{m k}\rangle,\qquad 
B^{(m)}_{l \nu} =
\langle w_{m l}, w_{m \nu} \rangle, 
\\
\label{vecchc} c^{(m)}_l & =& \langle w_{m l},
Qf \rangle,\qquad 
\hat{c}^{(m)}_l = \langle
w_{m l}, y \rangle, 
\\
\label{veccrz} r^{(m)}_l & =& \langle w_{m l}, Q
R_m \rangle,\qquad 
z^{(m)}_k =
a_{mk}, 
\\
\label{vechhh} h^{(m)}_\nu& =& a_{m \nu},\qquad
\hat{h} {}^{(m)}_\nu= \hat{a}_{m \nu} =
\langle y, t_{m \nu}\rangle, 
\end{eqnarray}
%
where $k,l \in K_{0m}$, $\nu\in K_{0m}^c$ and $\hat{a}_{mk}$ with $k,l
\in K_{0m}$, are defined in (\ref{coefest}). Then, one can rewrite an
exact system of linear equations (\ref{eqscprod}) as
$\mathbf{c}^{(m)}= \mathbf{A}^{(m)}\mathbf{z}^{(m)}+
\mathbf{B}^{(m)}\mathbf{h}^{(m)}+ \mathbf{r}^{(m)}$
and obtain its approximate version
%
%
\begin{equation}
\label{eqapproxsys} \hat{\mathbf{c}} {}^{(m)}= \mathbf{A}^{(m)}
\hat{\mathbf{z}} {}^{(m)}+ \mathbf{B}^{(m)} \hat{\mathbf{h}}
{}^{(m)}.
\end{equation}
Since matrix $\mathbf{A}^{(m)}$ is a nonnegative definite matrix of a
finite size, in order to guarantee that it is nonsingular, it is
sufficient to impose the following almost trivial assumption:
\begin{longlist}[(A5)]
\item[(A5)] Functions $w_{m k}= Q \varphi_{mk}$, $k \in K_{0m}$ are
    linearly independent.
\end{longlist}
Under assumption (A5), one has
%
%
\begin{eqnarray}
\label{eqsolutions} \mathbf{z}^{(m)}& = & \bigl(\mathbf{A}^{(m)}
\bigr)^{-1} \bigl(\mathbf{c}^{(m)}- \mathbf{B}^{(m)}
\mathbf{h}^{(m)}- \mathbf{r}^{(m)} \bigr),
\nonumber
\\[-8pt]
\\[-8pt]
\hat{\mathbf{z}} {}^{(m)}& = & \bigl(\mathbf{A}^{(m)}
\bigr)^{-1} \bigl(\hat{\mathbf{c}} {}^{(m)}- \mathbf
{B}^{(m)}\hat{\mathbf{h}} {}^{(m)} \bigr).\nonumber
\end{eqnarray}
Finally, for a given $m$, we set $\hat{a}_{mk}= \hat{z}{}^{(m)}_{k}$,
$k \in K_{0m}$ and estimate $f_{0,m}$ by the following wavelet linear
estimator
%
%
\begin{equation}
\label{hfom} \hat{f}_{0,m}(x) = \sum_{k \in K_{0m}}
\hat{a}_{mk}\varphi_{mk}(x),\qquad x \in[0,1].
\end{equation}

%
%
\begin{remark}[(Relation to nonparametric regression estimation based on spatially inhomogeneous data)]\label{relation}
We need to touch upon the relationship between the present paper and
the paper by \citet{anton} which considered nonparametric
regression estimation based on irregularly spaced data, in particular,
in the case when design density has zeros. The latter problem is the
well-known formulation which has been already studied extensively [see,
e.g., \citeauthor{gaiffas} (\citeyear{gaiffas,gaiffas1,gaiff2,gaiff3})]
and, indeed, can be considered as a trivial case of deconvolution with
the spatially inhomogeneous design studied in Section~\ref{secirreg}
with $Q$ being an identity operator. For this reason, the hybrid
estimator was proposed in \citet{anton}. However, due to the fact
that, in the regression set up, one observes function $f$ directly,
construction of the hybrid estimator is much more involved in the case
of an inverse problem than in the case of nonparametric regression. In
addition, the present paper provides the implementation of the hybrid
estimator and studies its performance via simulations which has never
been done previously since \citet{anton} considered only
theoretical construction of the hybrid estimator.
%
%
\end{remark}


\section{The risks of the estimators of the singularity-free and the
singularity-affected parts} \label{secriskfcf0}

In this section, we shall provide asymptotic expressions for the risks
of estimators (\ref{hfcmfixed}) and (\ref{hfom}) when resolution level
$m$ is a fixed, nonrandom quantity possibly dependent on $\varepsilon
\dvtx m = m(\varepsilon)$.


Let us first construct an asymptotic upper bound for the
singularity-free portion~(\ref{hfcmfixed}) of the estimator. Denote
%
%
\begin{equation}
\label{eqlamm} \lambda_{m}^{-2} = \sum
_{k \in K_{0m}^c} \lambda_{m k}^{-2}
\end{equation}
and observe that, under condition (\ref{eqlambdajk}), there exist
positive constants $C_{\lambda\alpha0}$ and $C_{\lambda\alpha}$
independent of $m$ such that
%
%
\begin{equation}
\label{eqlammexpr} C_{\lambda\alpha0} 2^{m(\beta+ \max(1, \alpha))} m^{\mathbb
{I}(\alpha
=1)} \leq
\lambda_{m}^{-2} \leq C_{\lambda\alpha} 2^{m(\beta+ \max(1, \alpha))}
m^{\mathbb
{I}(\alpha=1)}.
\end{equation}



%
%
\begin{lemma} \label{lemsingulfree}
Let $1 \leq p,q \leq\infty$, $s \geq\max{(1/2, 1/p)}$, ${s^*}= \min(s,
s')$ and assumptions \textup{(A1)--(A4)} hold. Let $\hat{f}_{c,m}$ be
given by (\ref{hfcmfixed}) where the nonrandom quantity $m =
m(\varepsilon)$ is such that $m_1 \leq m \leq J-1$, with $J$ defined in
(\ref{Jvalueo}). If
%
%
\begin{equation}
\label{eqmutau} \tau^2 = 4 C_u C_\lambda(
\sqrt{2} \chi+1)^2\quad\mbox{and}\quad \chi^2 \geq
\frac{4(\beta+ \max(1, \alpha))}{2 +\alpha+\beta},
\end{equation}
where $C_u$ and $C_\lambda$ are defined in (\ref{eqill-posedness}) and
(\ref{eqlambdajk}), respectively, then, as $\varepsilon\to0$,
%
%
\begin{equation}
\label{riskhfc} \sup_{f \in B_{p,q}^s (A)} \mathbb{E}\| \hat{f}_{c,m}-
f_{c,m}\|^2 \leq C \bigl( \lambda_{m}
^{-2} \varepsilon+ \Delta(\varepsilon) \bigl[\ln \bigl(\varepsilon
^{-1} \bigr) \bigr]^\rho \bigr).
\end{equation}
Here $\Delta(\varepsilon)$ is defined in (\ref{eqDDeps}),
%
%
\begin{equation}
\label{eqrho} \rho= \cases{ \displaystyle\mathbb{I} \biggl(\frac{2 s}{1-2/p} =
\frac{\beta+1}{\alpha-1} \biggr)+ \frac
{(1-\alpha)(2 - p)}{2 - \alpha p} \mathbb{I}(\alpha<1)
\mathbb{I}(p<2),
\vspace*{2pt}\cr
\qquad\mbox{if }\alpha\neq1,
\vspace*{4pt}\cr
\bigl(2 s^* + \beta+1 \bigr)^{-1} 2s^*,
\vspace*{2pt}\cr
\qquad\mbox{if }\alpha= 1}
\end{equation}
and $\mathbb{I}$ is the indicator function. Moreover, as $\varepsilon
\to0$,
%
%
\begin{equation}
\label{riskhfc4} \sup_{f \in B_{p,q}^s (A)} \mathbb{E}\|
\hat{f}_{c,m}- f_{c,m}\|^4 = o \bigl(\varepsilon
^{-2} \bigr).
\end{equation}
\end{lemma}


Now, we find upper bounds for the singularity-affected portion of the
estimator $\hat{f}_{0,m}(x)$. Recall that $w_{m k}= Q \varphi_{mk}$ and
let $\rho_m$ be such that
%
%
\begin{equation}
\label{eqwmk} C_{w1} \rho_m \leq\|w_{m k}\|
\leq C_{w2} \rho_m\qquad\mbox{if } k \in K_{0m}.
\end{equation}
Note that, since the set $K_{0m}$ contains at most $(2 D_0)$ indices,
$\rho_m$ satisfying condition (\ref{eqwmk}) can always be found. The
advantage of using the system of equations (\ref{eqapproxsys}) rests
upon the fact that matrix $\mathbf{A}^{(m)}$ is a finite-dimensional
positive definite matrix with all eigenvalues of order $\rho_m^2$. In
particular,
%
%
\begin{equation}
\label{eqAmrel} \bigl\|\mathbf{A}^{(m)}\bigr\| \leq C_{A1}
\rho_m^2,\qquad \bigl\| \bigl(\mathbf{A}^{(m)}
\bigr)^{-1} \bigr\| \leq C_{A2} \rho_m^{-2}
\end{equation}
for some positive constants $C_{A1}$ and $C_{A2}$ independent of $m$,
as it is shown in the proof of the following lemma which states the
rate of convergence of the singularity-affected portion of the
estimator.

%
%
\begin{lemma} \label{lemsingulaffect}
Let $1 \leq p,q \leq\infty$ and $s \geq\max(1/2, 1/p)$. Let assumptions
\mbox{\textup{(A1)--(A5)}} hold and there exists $C_{\rho\lambda}$, $0
< C_{\rho\lambda}< \infty$, independent of $m$, such that
%
%
\begin{equation}
\label{eqCrol} \rho_m^2 \geq C_{\rho\lambda}
\lambda_{m}^2,
\end{equation}
where $\lambda_{m}$ and $\rho_m$ are defined in (\ref{eqlamm}) and
(\ref{eqwmk}), respectively. Let $\alpha\geq1$ and also
%
%
\begin{eqnarray}
\label{eqlem21}\rho_m^{-4} \max_{l \in K_{0m}}
\sum_{k \in K_{0m}^c
}\lambda_{m k}^{-2}
\langle w_{m l}, w_{m k}\rangle^2 &\leq&
K_1 \lambda_{m}^{-2},
\\
\label{eqlem22}\max_{l \in K_{0m}} \sum_{j=m}^\infty
\sum_{k=0}^{2^j -1} \bigl|\langle w_{m l},
v_{j k}\rangle \bigr| &\leq& K_2 \rho_m^{2}
\end{eqnarray}
for some absolute constants $K_1$ and $K_2$ independent of $m$. Let
estimator $\hat{f}_{0,m}$ of~$f_{0,m}$ be given by (\ref{hfom}). Then,
for any $m$, $m_1 \leq m \leq J-1$, and some constant $C$ independent
of $m$ and $\varepsilon$, as $\varepsilon\to0$, one has
%
%
\begin{eqnarray}
\label{riskhfom} \sup_{f \in B_{p,q}^s (A)} \mathbb{E}\|
\hat{f}_{0,m}- f_{0,m}\|^2 &\leq& C \bigl(
\varepsilon\lambda_{m}^{-2} + 2^{-2m s'} \bigr),
\\
\label{riskhfom4} \sup_{f \in B_{p,q}^s (A)} \mathbb{E}\|
\hat{f}_{0,m}- f_{0,m}\|^4 &=& o \bigl(
\varepsilon ^{-2} \bigr).
\end{eqnarray}
\end{lemma}

Note that in order $\hat{f}_{m}= \hat{f}_{c,m}+ \hat{f}_{0,m}$
estimates $f$, one needs to start the estimator $\hat{f}_{c,m}$ in
(\ref{hfcmfixed}) at exactly the same resolution level at which the
linear estimator $\hat{f}_{0,m}$ in (\ref{hfom}) is constructed. Thus,
the choice of the lowest resolution level in (\ref{hfcmfixed}) is
driven by the choice of $m$ in (\ref{hfom}).

Let $m_0$ be such that
%
%
\begin{equation}
\label{eqmo} 2^{m_0} = \bigl(\varepsilon \bigl[\ln \bigl(
\varepsilon^{-1} \bigr) \bigr]^{\mathbb{I}(\alpha=1)} \bigr)^{-1/(2s' + \alpha+ \beta)},
\end{equation}
so that, for $\alpha\geq1$, one has
%
%
\begin{equation}
\label{eqtwoterms} \varepsilon\lambda_{m_0}^{-2} +
2^{-2m_0s'} \leq C \Delta(\varepsilon) \bigl[\ln \bigl(
\varepsilon^{-1} \bigr) \bigr]^\rho.
\end{equation}
%

The following statement delivers the total squared risk of the
estimator (\ref{totalest}) of $f$ if $D D_0>0$.

%
%
\begin{theorem} \label{thrisknonadapt}
Let $1 \leq p,q \leq\infty$ and $s \geq\max(1/2, 1/p)$. Let conditions
\mbox{(\ref{eqCrol})--(\ref{eqlem22})} and assumptions
\textup{(A1)--(A5)} hold with $\alpha\geq1$ and $D D_0 >0$ in
assumptions \textup{(A2)} and \textup{(A4)}. Consider estimator
(\ref{totalest}) of $f$ where $\hat{f}_{c,m}$ and $\hat{f}_{0,m}$ are
given by formulae (\ref{hfcmfixed}) and (\ref{hfom}), respectively. Let
$m=m_0$ where $m_0$ is defined in (\ref{eqmo}), $J$ be defined
in~(\ref{Jvalueo}) and let positive constants $\tau$ and $\chi$ satisfy
condition~(\ref{eqmutau}). Then, as $\varepsilon\to0$, one has
%
%
\begin{equation}
\label{eqrisknonadapt} \sup_{f \in B_{p,q}^s (A)} \mathbb{E}\|
\hat{f}_{m_0} - f \|^2 \leq C \Delta(\varepsilon) \bigl[\ln
\bigl(\varepsilon^{-1} \bigr) \bigr]^\rho,
\end{equation}
where $\Delta(\varepsilon)$ and $\rho$ are defined in (\ref{eqDDeps})
and (\ref{eqrho}), respectively.
\end{theorem}

Validity of Theorem \ref{thrisknonadapt} follows directly from Lemmas
\ref{lemsingulfree} and \ref{lemsingulaffect} and
inequality~(\ref{eqtwoterms}).

Note that the value of $m_0$ depends on the unknown parameters $s$ and
$p$ of the Besov space, hence, in general, estimator $\hat{f}_{m_0}$ in
(\ref{eqrisknonadapt}) is not adaptive. However, if $D=D_0 =0$, then
$f_{0,m}= \hat{f}_{0,m}\equiv0$, $\hat{f}= \hat{f}_{c,m}$ and one can
choose $m=m_1$ in $\hat{f}_{c,m}$ using formula (\ref{eqm1val}), so
that the estimator is adaptive. In this case, convergence rates of
$\hat{f}$ are given entirely by Lemma \ref{lemsingulfree}. In
particular, the following corollary is valid.

%
%
\begin{corollary} \label{corsingfree}
Let $D=D_0 =0$ and assumptions of Lemma \ref{lemsingulfree} hold.
Consider estimator $\hat{f}_m = \hat{f}_{c,m}$ given by
(\ref{hfcmfixed}) with $m=m_1$. Then
%
%
\begin{equation}
\label{eqrtotalsingfree} \sup_{f \in B_{p,q}^s (A)} \mathbb{E}\|
\hat{f}_{m_1} -f \|^2 \leq C \Delta(\varepsilon) \bigl[\ln
\bigl(\varepsilon^{-1} \bigr) \bigr]^\rho\qquad (\varepsilon
\to0),
\end{equation}
where $\Delta(\varepsilon)$ and $\rho$ are defined in (\ref{eqDDeps})
and (\ref{eqrho}), respectively.
\end{corollary}

If $D D_0 >0$, then it is necessary to construct an adaptive estimator
of $f$. Note that this is not an easy task. Expanding the system of
equations in (\ref{eqscprod}) so that it includes not only the scaling
but also the wavelet coefficients will compromise uniformity of
eigenvalues of matrix $\mathbf{A}^{(m)}$ [see (\ref{eqAmrel})] which
are ensured by positive-definiteness and finite size of
$\mathbf{A}^{(m)}$. On the other hand, introducing a penalty on the
solution does not help either since, for any $m$, the system involves
the unknown bias term $\mathbf{r}^{(m)}$. For this reason, in order to
choose parameter $m$, we apply Lepski's method since it allows us to
eliminate the bias inherent to the system of
equations~(\ref{eqapproxsys}).

\section{Adaptive estimation in the presence of singularity}
\label{secadaptive}

In order to construct an adaptive estimator of $f$ in the presence of
singularity, we shall use the technique of optimal tuning parameter
selection pioneered by \citeauthor{le1990} (\citeyear{le1990,le1991})
and further exploited in \citet{lepskii2} and
\citet{lepskii1}. The idea behind this technique is to construct
estimators for various values of the tuning parameter in question ($m$,
in our case), and then choose an optimal value of the tuning parameter
by regulating the differences between the estimators constructed with
different values of the parameter.

In particular, for various values of $m$, we construct versions of the
system of equations, obtain\vspace*{-1pt} values of $\hat{\mathbf{z}}{}^{(m)}$ in
(\ref{eqsolutions}) and use them as $\hat{a}_{mk} =
\hat{z}{}^{(m)}_{k}$, $k \in K_{0m}$, in~(\ref{hfom}). Then, for
various values of $m$, we obtain estimators $\hat{f}_{m}$ of $f$ using
formula (\ref{totalest}) where $\hat{f}_{0,m}$ and $\hat{f}_{c,m}$ are
of the forms (\ref{hfom}) and (\ref{hfcmfixed}), respectively, and $m$
is the lowest resolution level of $\hat{f}_{c,m}$. After that, we
choose the ``best possible'' resolution level $\hat{m}$ and consider
estimator $\hat{f}_{\hat{m}}$ as the final estimator. The choice of the
resolution level $\hat{m}$ is driven by the singularity-affected
portion of $f$ rather than the zero-free portion as it is described
below.

For any resolution level $m >0$, we define a neighborhood $\Omega_m$ of
$x_0$ as
%
%
\begin{eqnarray}\label{eqOmneighborhood}
\Omega_m &=& \bigl\{x\dvtx \min({L_\varphi}-D_0,{L_\psi}-D) < 2^m
(x-x_0)
\nonumber\\[-8pt]\\[-8pt]
&&\hspace*{67pt} < \max({U_\varphi}+D_0,{U_\psi}+D)\bigr\},\nonumber
\end{eqnarray}
where $\operatorname{supp} \varphi= ({L_\varphi}, {U_\varphi})$ and
$\operatorname{supp} \psi= ({L_\psi}, {U_\psi})$. Observe that
$\Omega_m$ is designed so that $\operatorname{supp}(f_{0,m})
\subseteq\Omega_m$, $\operatorname{supp}(\hat{f}_{0,m})
\subseteq\Omega_m$ and $\Omega_j \subset\Omega_m$ if $j > m$.

Choose $m = \hat{m}$ such that $m_1 \leq m \leq J-1$, where $J$ is
defined in (\ref{Jvalueo}) and
%
%
\begin{eqnarray}\label{mopt}
\hat{m}&=& \min \bigl\{m\dvtx\bigl\| (\hat{f}_{m}-
\hat{f}_{j}) \mathbb{I}(\Omega_m) \bigr\|^2 \leq
\kappa^2 \varepsilon\ln \bigl(\varepsilon^{-1} \bigr)
\lambda_j^{-2}
\nonumber\\[-8pt]\\[-8pt]
&&\hspace*{110pt} \mbox{ for all } j \in[m, J-1] \bigr\},\nonumber
\end{eqnarray}
where $\kappa>0$ is a constant which is defined below.

The construction of $\hat{m}$ is based on the following idea. Note that
when \mbox{$m = \hat{m}\leq m_0$,} one has
%
%
\begin{equation}
\label{Lepskii-expl} \mathbb{E} \| \hat{f}_{\hat{m}}- f \|^2
\leq2 \bigl[\mathbb{E}\| \hat{f}_{\hat{m}}- \hat{f}_{m_0}
\|^2 + \mathbb{E}\| \hat{f}_{m_0}- f \|^2
\bigr].
\end{equation}
The first component in (\ref{Lepskii-expl}) is small due to the
definition of the resolution level $\hat{m}$ while the second component
is calculated at the optimal resolution level $m_0$ and, hence, tends
to zero at the optimal convergence rate (up to a logarithmic factor).
On the other hand, if $m=\hat{m}> m_0$, then there should exist $j > m$
such that $\| (\hat{f}_{m}- \hat{f}_{j}) \mathbb{I}(\Omega_m) \|^2 >
\kappa^2 \varepsilon\ln(\varepsilon^{-1})\lambda_j^{-2}. $ The
following Lemma shows that the probability of this event is
infinitesimally small provided $\kappa$ is large enough.



%
%
\begin{lemma} \label{lemmoptdev}
Let $m_0$ and $\hat{m}$ be given by expressions (\ref{eqmo}) and
(\ref{mopt}) 
and the conditions of Theorem \ref{thrisknonadapt} hold. Let $C_\kappa=
2^{11} D_0^2 C_{A2}^2 \max(C_{\rho\lambda}^{-1} C_{w2}^2, C_t K_1)$
where\vspace*{1pt} constants $C_{A2}$, $C_{\rho\lambda}$, $C_{w2}$, $C_t$ and $K_1$
are defined in (\ref{eqAmrel}), (\ref{eqCrol}), (\ref{eqwmk}),
(\ref{eqtjk}) and~(\ref{eqlem21}), respectively. If
%
%
\begin{equation}
\label{eqkappamucond} \kappa\geq\max \bigl(d^2 C_\kappa,
2^6 \bigl(d^2 +1 \bigr) C_t \bigr),\qquad
\chi^2 \geq d^2 + 2/(2 + \alpha+ \beta),
\end{equation}
then, as $\varepsilon\to0$,
%
%
\begin{equation}
\label{largedevm} \mathbb{P} (\hat{m}> m_0 )\leq C
\varepsilon^{d^2}.
\end{equation}
\end{lemma}


Lemma \ref{lemmoptdev} confirms that indeed $m= \hat{m}$ can be chosen
as the lowest resolution level in the nonlinear portion of the
estimator, so that we estimate $f$ by
%
%
\begin{equation}
\label{eqnonlin} \hat{f}(x) = \hat{f}_{0,\hat{m}} (x) +
\hat{f}_{c,\hat{m}}(x), \qquad x \in[0,1],
\end{equation}
where $\hat{f}_{0,m}(x)$ and $\hat{f}_{c,m}(x)$ are defined in
(\ref{hfom}) and (\ref{hfcmfixed}), respectively. The following
statement confirms that the wavelet nonlinear estimator $\hat{f}$ given
by (\ref{eqnonlin}) indeed attains (up to a logarithmic factor) the
asymptotic minimax lower bounds obtained in Theorem \ref{thlower}.


%
%
\begin{theorem} \label{thriskadapt}
Let $1 \leq p,q \leq\infty$ and $s \geq\max(1/2, 1/p)$. Let conditions
\mbox{(\ref{eqCrol})--(\ref{eqlem22})} and assumptions \textup{(A1)--(A5)}
hold with $\alpha\geq1$ and $D D_0 >0$ in assumptions \textup{(A2)} and
\textup{(A4)}. Consider the estimator (\ref{totalest}) of $f$ where
$\hat{f}_{c,m}$ and $\hat {f}_{0,m}$ are given by formulae
(\ref{hfcmfixed}) and (\ref{hfom}), respectively. Let $m=\hat{m}$ where
$\hat{m}$ is defined in~(\ref{mopt}). Let $J$ be defined in
(\ref{Jvalueo}) and $\kappa$ and $\chi$ be such that
%
%
\begin{equation}
\label{eqcondth3} \kappa\geq\max(4 C_\kappa, 320 C_t ),
\qquad \chi^2 \geq4 + 2/(2 + \alpha+ \beta),
\end{equation}
where $C_\kappa$ and $C_t$ are defined in Lemma \ref{lemmoptdev} and
formula (\ref{eqtjk}), respectively.
Then, 
%
%
\begin{equation}
\label{eqriskadapt} \sup_{f \in B_{p,q}^s (A)} \mathbb{E}\|
\hat{f}_{\hat{m}} - f \|^2 \leq C \Delta(\varepsilon) \bigl[\ln
\bigl(\varepsilon^{-1} \bigr) \bigr]^{1 + \mathbb{I}(\alpha=1)}\qquad (\varepsilon
\to0),
\end{equation}
where $\Delta(\varepsilon)$ is defined in (\ref{eqDDeps}).
\end{theorem}

%
%
\begin{remark}[(Logarithmic factor in convergence rates)]\label{nonadapt}
Note that in (\ref{eqrho}), $\rho\neq0$ only if $2 s (\alpha-1) =
(\beta+1)(1-2/p)$, or $\alpha= 1$, or $\alpha<1$ and $p \leq2$. The
latter shows that the lower bounds for the risk in
Theorem~\ref{thlower} cannot be made tighter, at least, in the case
when $\alpha>1$. Theorems~\ref{thlower} and \ref{thrisknonadapt} and
Corollary~\ref{corsingfree} demonstrate that estimator (\ref{totalest})
attains the asymptotically optimal convergence rates if $2s(\alpha-1)
\neq(\beta+1)(1-2/p)$ and $\alpha>1$, or if $\alpha<1$ and $p>2$.
Otherwise, estimator (\ref{totalest}) is asymptotically near-optimal up
to a logarithmic factor.
%
However, in Theorem~\ref{thriskadapt} the risk of the adaptive
estimator is always within a logarithmic factor $ [\ln(\varepsilon
^{-1})]^{1 + \mathbb{I} (\alpha=1)}$ of the minimax risk. The latter is
due to application of Lepski method. Note that in spite of the fact
that we are using the integrated mean squared error, Lepski method is
applied locally and, hence, leads to an extra log-factor in the risk,
as it usually happens with application of Lepski method to pointwise
estimation.
\end{remark}


\section{Examples}
\label{secexamples}

\subsection{Deconvolution with a spatially inhomogeneous kernel}
\label{secdeconv}

Consider problem~(\ref{eqmaineq}) with operator $Q$ of the form
%
%
\begin{equation}
\label{eqex1} (Q f) (x) = \mu(x) \int_0^1
q(x-t) f(t) \,dt,\qquad x \in[0,1],
\end{equation}
where functions $\mu(x)$, $q(x)$ and $f(x)$ are periodic and both
$q(x)$ and $\mu(x)$ are completely known. Problem (\ref{eqmaineq}) is
equivalent to the following statistical problem:
%
%
\begin{equation}
\label{eqheteroscedastic} y_i = \mu(i/n) \int_0^1
q(i/n -t) f(t)\,dt + \sigma\xi_i,\qquad i=1,\ldots, n,
\end{equation}
where $\xi_i$ is a white Gaussian noise and $\varepsilon= \sigma^2/n$.
Equation of the form (\ref{eqex1}) can appear when one observes a
convolution $Y(x)$ of the known kernel $q$ with the unknown function of
interest $f$ and a known heteroscedastic noise $ \sqrt{\varepsilon}
\gamma(x) W(x)$, so that $\mu(x) = [\gamma(x)]^{-1}$ and $Y(x) =
\gamma(x) y(x) = y(x)/\mu(x)$. In this case, equation~(\ref{eqheteroscedastic}) takes the form
%
%
\begin{equation}
\label{eqdifnoise} Y_i = \int_0^1
q(i/n -t) f(t)\,dt + \sigma\gamma(i/n) \xi_i,\qquad i=1,\ldots, n.
\end{equation}


If $\mu(x)$ is uniformly bounded above and below, in principle, spatial
inhomogeneity of operator $Q$ in (\ref{eqex1}) can be ignored. Below we
consider the case when the former is not true since $\mu(x)$ vanishes
at some point $x_0\in(0,1)$, in particular,
%
%
\begin{equation}
\label{eqmu} C_{\mu1} |x-x_0|^{\alpha} \leq
\mu^2 (x) \leq C_{\mu
2} |x-x_0|^{\alpha}
\end{equation}
for some for some positive constants $\alpha$, $C_{\mu1}$ and
$C_{\mu2}$ independent of $x_0$ and $x$. Therefore, the version of the
problem studied in the present paper can be described as locally
extreme noise which occurs when the degree of spatial inhomogeneity is
high.
Direct calculations show that
\[
\bigl(Q^* h \bigr) (z) = \int_0^1 q(x-z) h(x)
\mu(x) \,dx.
\]
Hence, functions $u_{j k}$ and $t_{m k}$ are solutions of the equations
\[
Q^* u_{j k}= \psi_{jk},\qquad Q^* t_{m k}=
\varphi_{mk}.
\]
%
Consider $\omega\in\mathbb{Z}$ and let $e_\omega(t)= e^{i 2\pi\omega
t}$, $t \in[0,1]$, be the periodic Fourier basis on $[0,1]$. Denote
$U_{jk}(x) = u_{j k}(x) \mu(x)$ and, similarly, $T_{mk}(x) = t_{m k}(x)
\mu(x)$ and introduce Fourier coefficients $q_\omega= \langle q,
e_\omega\rangle$, $U_{jk\omega}= \langle U_{jk}, e_\omega\rangle$,
$T_{mk\omega}= \langle T_{mk}, e_\omega\rangle$, $\psi_{jk\omega}=
\langle\psi_{jk}, e_\omega\rangle$ and $\varphi_{mk\omega}=
\langle\varphi_{mk}, e_\omega\rangle$.
Then, $U_{jk\omega}= [\bar{q}_{\omega}]^{-1} \psi_{jk\omega}$,
$T_{mk\omega}= [\bar{q}_{\omega} ]^{-1} \varphi_{mk\omega}$ and
%
%
\begin{eqnarray}\label{eqUjkTmk}
U_{jk}(x) &=& \sum_{\omega\in\mathbb{Z}}[
\bar{q}_{\omega}]^{-1} \psi_{jk\omega}e_\omega(x),
\nonumber\\[-8pt]\\[-8pt]
T_{mk}(x) &=& \sum_{\omega\in
\mathbb{Z}}[\bar{q}_{\omega}]^{-1} \varphi_{mk\omega}e_\omega(x),\nonumber
\end{eqnarray}
where $\bar{q}_{\omega}$ is the complex conjugate of $q_\omega$.
Moreover, estimators $\hat{a}_{mk}$ and $\hat{b}_{jk}$ in~(\ref{coefest}) can be constructed using Fourier wavelet transform
suggested in \citet{jpkr}. Indeed, if $Y(x) =
y(x)/\mu(x)$ and $Y_\omega$ are Fourier coefficients of function
$Y(x)$, then
%
%
\begin{eqnarray}\label{eqwavcoefex1}
\hat{b}_{jk} &=& \langle y, u_{j k}\rangle= \langle Y, U_{jk}\rangle= \sum_{\omega
\in\mathbb{Z}}[\bar{q}_{\omega}]^{-1} \psi_{jk\omega}
\overline{Y_\omega},
\nonumber\\[-8pt]\\[-8pt]
\hat{a}_{mk}&=& \sum_{\omega\in
\mathbb{Z}}[ \bar{q}_{\omega} ]^{-1}
\varphi_{mk\omega} \overline{Y_\omega}.\nonumber
\end{eqnarray}
In the case of the statistical experiment described in formula (\ref
{eqheteroscedastic}), Fourier coefficients are replaced by the discrete
Fourier transform.

Note that application of the wavelet-vaguellete methodology to
deconvolution with heteroscedastic noise (\ref{eqdifnoise}) appears
very reasonable. Really, if noise level $\gamma(x)$ is such that $\int
\gamma^2 (x) \,dx < \infty$, then formula (\ref{eqwavcoefex1}) implies
that wavelet coefficients are estimated using Fourier transform of the
measured signal $Y$ in (\ref{eqdifnoise}) and then thresholded taking
into account the local noise level. Indeed, it is easy to observe that
$\lambda_{j k}$ in (\ref{eqlamjkex1}) is\vspace*{-1pt} such that $\lambda_{j k}^{-2}
\asymp\| U_{jk} \|^2 \int\gamma^2(x) \mathbb{I}(x
\in\operatorname{supp} U_{jk}) \,dx$. If $\int\gamma^2 (x) \,dx = \infty$
in the vicinity of some point $x_0$, the natural strategy suggested
above cease working and one needs another means for estimating scaling
and wavelet coefficients in the neighborhood of $x_0$. In this
situation one has to apply the hybrid estimator constructed in
Sections~\ref{secestimationstrategies} and \ref{secadaptive}. In
Section~\ref{secsimulations}, we provide a detailed description of the
computational algorithm for this task.

If function $\mu(x)$ is not completely known and is estimated from
data, matrices $\mathbf{A}^{(m)}$ and $\mathbf{B}^{(m)}$ as well as
vector $\hat{\mathbf{c}}{}^{(m)}$ will be subjected to additional
errors which have to be accounted for by using regularization
techniques designed for the inverse problems with errors in the
operator [see, e.g., \citet{engl} and \citet{hoffmann}].

In order to find $\lambda_{j k}\asymp\| v_{j k}\|$ in assumption (A1)
and verify assumptions \mbox{(A1)--(A5)}, we impose the following conditions
on the kernel $q$ and mother and father wavelets $\psi$ and $\varphi$.
\begin{longlist}[(E3)]
\item[(E1)] Kernel $q(x)$ is $(r-2)$ times continuously
differentiable on $[0,1]$ and $r_1 >r\geq1$ times differentiable
outside the neighborhood of jump discontinuities of $q^{(r-1)}$
with $q^{(r)}$ and $q^{(r_1)}$ uniformly bounded. The value $r=1$
corresponds to the case when $q$ itself has jump discontinuities.

\item[(E2)] Fourier coefficients $q_\omega$ of $q$ are such that
$C_{q1} (|\omega|+1)^{-r} \leq|q_\omega| \leq C_{q2}
(|\omega|+1)^{-r} $ for some positive constants $C_{q1}$ and
$C_{q2}$ independent of $\omega$.

\item[(E3)] Let $\psi$ be $r_0$-regular, $r_0$ times continuously
differentiable wavelet function with the bounded support, where
$r_0> \max(r, r_1)$.

\item[(E4)] Kernel $q$ is such that functions $U_{jk}$ and $T_{jk}$
defined in (\ref{eqUjkTmk}) have bounded supports of the lengths
proportional to $2^{-j}$ and centered at $2^{-j}
k\dvtx\break\operatorname {supp}(U_{jk}) = (2^{-j} (k-d_U), 2^{-j}
(k+d_U))$ and $\operatorname{supp}(T_{jk}) = (2^{-j} (k-d_T),
2^{-j} (k+d_T))$.
\end{longlist}

There are many functions $q$ satisfying conditions above, among them,
for example,
%
%
\begin{eqnarray}\label{eqvariousq}
q_1(x) &=& \sum_{k \in\mathbb{Z}}
\exp\bigl(-\lambda|x+k|\bigr)\quad\mbox{and}
\nonumber\\[-8pt]\\[-8pt]
q_2(x) &=& \sum_{k \geq0} \exp \bigl(-\lambda(x+k) \bigr) (x+k)^N\nonumber
\end{eqnarray}
with $r=2$ for $q_1(x)$ and $r=N+1$ for $q_2(x)$. Note that under
assumptions \mbox{(E1)--(E4)}, one has $\|u_{j k}\|^2 = \infty$
if $|k - k_{0j}| < d_U$ and $\alpha\geq1$ and
%
%
\begin{eqnarray}\label{eq:lamjkex1}
\|u_{j k}\|^2 & \asymp&
\mu^{-2} \bigl(2^{-j} k \bigr) \|U_{jk}
\|^2 \asymp\bigl| 2^{-j} k - x_0\bigr|^{-\alpha}
\sum_{\omega\in\mathbb{Z}}|\bar{q}_{\omega}|^{-2} |
\psi_{jk\omega}|^2
\nonumber
\\[-8pt]
\\[-8pt]
&\leq& C 2^{j \alpha} \bigl[|k - k_{0j}|^{\alpha} +1
\bigr]^{-1} \sum_{\omega\in\mathbb{Z}} \bigl(|
\omega|^{2r} + 1 \bigr) |\psi_{jk\omega}|^2,
\nonumber
\end{eqnarray}
otherwise. Due to conditions (E3) and (E4) and periodicity of
$\psi_{jk}$, using integration by parts $r$ times,
we obtain 
\begin{eqnarray*}
\psi_{jk\omega}&=& \int_0^1
\psi_{jk}(x) e^{i 2\pi\omega x} \,dx
\\
&=& 
2^{j(r+1/2)} (- 2\pi i \omega)^r \int
_0^1 \psi^{(r)} \bigl(2^j
x - k \bigr) e^{i 2\pi\omega x} \,dx,
\end{eqnarray*}
so that
$ 2^{-2jr} \sum_{\omega\in\mathbb{Z}}(|\omega|^{2r} + 1) |\psi
_{jk\omega}|^2 \leq C ( 2^{-2jr} \|\psi\|^2 + \|\psi^{(r)}\|^2). $
Therefore,
%
%
\begin{eqnarray}\label{equjknormex1}
\|u_{j k}\|^2 \leq C 2^{j(\alpha+2r)}
\bigl[|k - k_{0j}|^\alpha+1 \bigr]^{-1} \asymp
\lambda_{j k}^{-2},
\nonumber\\[-8pt]\\[-8pt]
\eqntext{\mbox{if }|k -k_{0j}|\geq
d_U \mbox{ or }  \alpha<1}
\end{eqnarray}
and $\|u_{j k}\| = \infty$ otherwise. Similar inequality can be proved
for $\| t_{m k}\|^2$.

The following proposition shows that, indeed, $\lambda_{j k}$ is
defined by expression (\ref{equjknormex1}) and that, under conditions
(E1)--(E4), assumptions (A1)--(A5) hold.

%
%
\begin{proposition} \label{propexample1}
Let $\mu(x)$ satisfy condition (\ref{eqmu}). Then, under assumptions
\textup{(E1)--(E4)} with $r_1$ such that $2 r_1 +1 > 2r + \alpha$, one has
%
%
\begin{equation}\label{eqlamjkex1}
\lambda_{j k}^2 \asymp2^{-j(2r+ \alpha)}
\bigl[|k - k_{0j} |^\alpha+1 \bigr].
\end{equation}
Furthermore, assumptions \textup{(A1)--(A5)} hold with $D=D_0=0$ if $\alpha<1$,
and with $D= d_U$ and $D_0 = d_T$ if $\alpha\geq1$. Conditions
(\ref{eqCrol})--(\ref{eqlem22}) of Lemma \ref{lemsingulaffect} are also
valid.
\end{proposition}

Due to Proposition \ref{propexample1}, all statements and constructions
in Sections~\ref{seclowbounds}--\ref{secadaptive} can be applied to
equation (\ref{eqmaineq}) with operator $Q$ given in (\ref{eqex1}). In
particular, Theorems~\mbox{\ref{thlower}--\ref{thriskadapt}} can be utilized
with $\beta= 2r$.
By direct comparison with, for example, \citet{jpkr}, one
can see that if $\alpha=0$, so that the problem is spatially
homogeneous, then the rates of convergence in Theorems
\ref{thlower}--\ref{thriskadapt} coincide with the usual convergence
rates exhibited in deconvolution problems with white noise.


\subsection{Deconvolution with spatially inhomogeneous design}
\label{secirreg}

Consider the problem of deconvolution when measurements are irregularly
spaced. In particular, let $g$ be a sampling p.d.f. with the corresponding\vadjust{\goodbreak}
c.d.f. $G$. Due to irregular design, operator $Q$ can be presented as
%
%
\begin{equation}
\label{eqex2} (Q f) (x) = \int_0^1 q
\bigl(G^{-1}(x)-t \bigr) f(t) \,dt,
\end{equation}
where $G^{-1}$ is the inverse of $G$.
In this case, equation (\ref{eqmaineq}) can be viewed as an idealized
version of the equation
%
%
\begin{equation}
\label{eqirregdesign} y_i = \int_0^1
q(x_i -t) f(t)\,dt + \sigma\xi_i,\qquad i=1,\ldots, n,
\end{equation}
where $\varepsilon= \sigma^2/n$, $\xi_i$ is a white Gaussian noise,
observation points $x_i$, $i=1,\ldots, n$, are such that $G(i/n) = x_i$
and $x_i$'s and $\xi_j$'s are independent. Then, \mbox{$G^{-1}(i/n)= x_i$},
and the right-hand side of (\ref{eqmaineq}) with operator $Q$ given by
(\ref{eqex2}) provides a continuous equivalent of the statistical
problem (\ref{eqirregdesign}).

In what follows, we assume that functions $q(x)$ and $f(x)$ are
periodic and both $q(x)$ and $g(x)$ are completely known. In this case,
the conjugate operator $Q^*$ is of the form
\[
\bigl(Q^* h \bigr) (z) = \int_0^1 q
\bigl(G^{-1}(x)-z \bigr) h(x) \,dx.
\]
It is pretty straightforward to show that $u_{j k}(x) =
U_{jk}(G^{-1}(x))$ and $t_{m k}(x) = T_{mk}(G^{-1}(x))$ where, as
before, $U_{jk}(\cdot)$ and $T_{mk}(\cdot)$ are given by formula
(\ref{eqUjkTmk}). Wavelet coefficients $\hat{b}_{jk}$ and
$\hat{a}_{mk}$ can also be estimated in a manner similar to Example~1.
Indeed, if $Y(x) = y(G(x))$ and $Y_\omega$ are Fourier coefficients of
$Y$, then $\hat{b}_{jk}$ and $\hat{a}_{mk}$ can be evaluated using
formula (\ref{eqwavcoefex1}).

If design density $g$ is unknown, then both $g(x)$ and $G(x)$ have to
be estimated from observations $x_i$, $i=1,\ldots, n$. The latter will
lead to additional errors in estimating wavelet coefficients
$\hat{b}_{jk}$ and $\hat{a}_{mk}$ as well as entries of matrices
$\mathbf{A}^{(m)}$ and $\mathbf{B}^{(m)}$ and vector
$\hat{\mathbf{c}}{}^{(m)}$. The issue of additional errors has to be
addressed by using, for example, regularization techniques [see, e.g.,
\citet{engl}  and \citet{hoffmann}].

We assume that design density $g(x)$ has a single zero of order
$\alpha$ at $x_0$, that is, $g(x_0 + x) |x|^{-\alpha} \to C_g$ as $x
\to0$. The latter implies that there exist some absolute constants
$C_{g1}$ and $ C_{g2}$ such that, for any $x$, one has
%
%
\begin{equation}
\label{eqgineq} C_{g1} |x-x_0|^{\alpha} \leq g(x)
\leq C_{g2} |x-x_0|^{\alpha}.
\end{equation}
Thus, we are considering the case of extremely inhomogeneous design
which can be described also as a local data loss.

In this case, under conditions (E1)--(E4), similarly to Example 1, one
has $\|u_{j k}\|^2 = \int g^{-1} (x) U_{jk}^2 (x) \,dx$. Hence,
identically to (\ref{equjknormex1}), one has
$\|u_{j k}\|^2 \leq C 2^{j (\alpha+ 2r)} |k - k_{0j}|^{-\alpha}$ if $|k
- k_{0j}|\geq d_U$ or $\alpha<1$, and $\|u_{j k}\|^2 = \infty$ if $|k -
k_{0j}| < d_U$ and $\alpha\geq1$.
Moreover, by simple modifications of the proof of
Proposition~\ref{propexample1}, it easy to show that the following
statement is valid.

%
%
\begin{proposition} \label{propexample2}
Let $g (x)$ satisfy condition (\ref{eqgineq}). Then, under assumptions
\textup{(E1)--(E4)} with $r_1$ such that $2 r_1 +1 > 2r + \alpha$, the value of
$\lambda_{j k}$ is given by formula~(\ref{eqlamjkex1}).
Furthermore, assumptions \textup{(A1)--(A5)} hold with $D=D_0=0$ if $\alpha<1$,
and with $D= d_U$ and $D_0 = d_T$ if $\alpha\geq1$. Conditions
(\ref{eqCrol})--(\ref{eqlem22}) of Lemma \ref{lemsingulaffect} are also
valid.
\end{proposition}

Again, analogously to Section~\ref{secdeconv}, due to Proposition
\ref{propexample2}, all statements and constructions in
Sections~\ref{seclowbounds}--\ref{secadaptive} can be applied to equation
(\ref{eqmaineq}) with operator $Q$ given in (\ref{eqex2}). In
particular, Theorems~\ref{thlower}--\ref{thriskadapt} can be used with
$\beta= 2r$.
If $q$ is the Dirac delta function, then $(Q f)(x) = f(G^{-1}(x))$,
then $r=0$ and the problem reduces to regression estimator based on
spatially inhomogeneous data studied in \citet{anton}. In this
case, the rates of convergence coincide with the minimax convergence
rates derived therein.

%
%
\begin{remark}[(Irregularly spaced observations and heterogeneous noise)]\label{equivalence}
It follows from examples in Sections~\ref{secdeconv} and \ref{secirreg}
that there is a direct correspondence between deconvolution with
irregularly spaced measurements and deconvolution with heterogeneous
noise. In particular, as far as convergence rates are concerned, the
squared noise level acts in a similar way to the inverse of the design
density and both are equivalent, in some way, to a multiplicative
factor in the convolution operator.
\end{remark}

\section{Simulation study and real data application}

\label{secsimulations}

\subsection{Simulation study}

In order to assess finite sample properties of the proposed methodology
and, in particular, performance of the hybrid estimator, we carried out
a small simulation study. We limited our attention to the deconvolution
in the presence of heteroscedastic noise described in
Section~\ref{secdeconv}. Specifically, we considered $q(x)=q_1 (x)$
with $\lambda= 5$ where $q_1 (x)$ is defined in (\ref{eqvariousq}). We
chose one of the standard test functions, \texttt{blip}, as the true
function $f(x)$. Function $\mu(x)$ in (\ref{eqex1}) is of the form
$\mu(x) = |h^{-1}(x - x_0)|^{\alpha/2} \mathbb{I}(|x - x_0| \leq h) +
\mathbb{I} (|x - x_0|
> h)$ with $x_0 = 1/3$ and $h=1/6$, so that condition (\ref{eqmu})
holds. We generated data using equation (\ref{eqdifnoise}) with
$\gamma(x) = \mu^{-1} (x)$ and $\sigma= 0.02$, in particular,
%
\[
Y_i = H(i/n) + \sigma\gamma(i/n) \xi_i,\qquad i=1,
\ldots, n
\]
with
\[
H(x) = \int_0^1 q(x
-t)f(t)\,dt.
\]
We evaluate noise intensity by the common in signal processing,
signal-to-noise ratio (SNR) which is defined as $\mathrm{SNR}= \sqrt{n}
\operatorname{std}(f)/(\|\gamma\|*\sigma)$ where $\| \gamma\|$ is the
$L^2$-norm of $\gamma$ and $\operatorname{std}(p)$ is the standard
deviation of $p$ for any function~$p(x)$.\vadjust{\goodbreak}
%
%
\begin{figure}[t]

\includegraphics{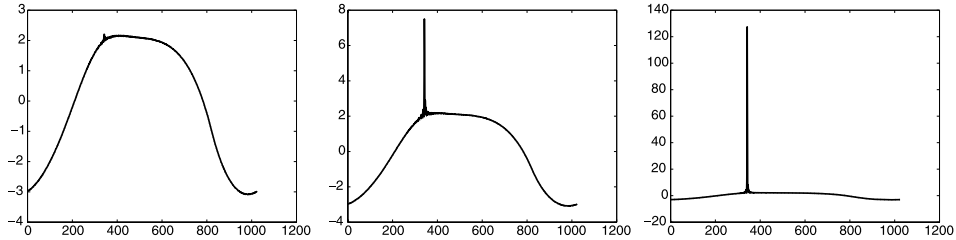}

\caption{True function $H$ (dotted line) and observed data (solid line)
with $\alpha= 1$ (left), $\alpha= 2$ (middle) and $\alpha=3$ (right).
Here, $\mathrm{SNR}= 0.8848$ for $\alpha= 1$, $\mathrm{SNR}= 0.0808$
for $\alpha=2$, $\mathrm{SNR}= 0.0183$ for $\alpha= 2.5$ and
$\mathrm{SNR}= 0.0040$ for $\alpha=3$.} \label{figdata}
\end{figure}
%
%
\begin{figure}[b]

\includegraphics{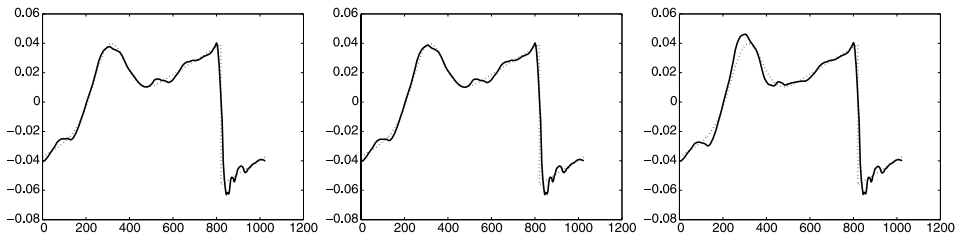}

\caption{Thresholded wavelet--vaguelette deconvolution estimator with
$n=1024$. True regression $f$ (dotted line) and its estimated value
(solid line) with $\alpha= 0$ (left), $\alpha= 1$ (middle) and $\alpha= 2$ (right).}\label{figwavelet}
\end{figure}

We used WaveLab package for Matlab and carried out simulations using
degree 8 Daubechies wavelets and $n=1024$. In order to obtain
estimators of wavelet and scaling coefficients, we generated wavelet
and scaling functions $\psi_{jk}$ and $\varphi_{mk}$ using
\texttt{MakeWavelet} command and obtained the respective matrix of the
Fourier coefficients. Subsequently, we found estimators of wavelet and
scaling coefficients using formula (\ref{eqwavcoefex1}) with $
Y_\omega$ being discrete Fourier transform of vector $Y$ in
(\ref{eqdifnoise}). We generated values of $\lambda_{j k}^{-2}$ using
equation (\ref{eq:lamjkex1}) and used them for hard thresholding. Due to
relatively small value of $n$, we did not use block thresholding
described in Section~\ref{secestimationstrategies}. By applying inverse
wavelet transform to the thresholded wavelet coefficients, we obtained
deconvolution estimator $\hat{f}$.

We evaluated performance of the estimators for $n=1024$ and different
values of~$\alpha$. As it is expected, when $\alpha$ is growing, the
SNRs are decreasing and the quality of observed data is declining.
Figure~\ref{figdata} demonstrates observed data for various values of
$\alpha$. The corresponding signal-to-noise ratios are $\mathrm{SNR}=
0.8848$ for $\alpha= 1$, $\mathrm{SNR}= 0.0808$ for $\alpha=2$ and
$\mathrm{SNR}= 0.0040$ for $\alpha=3$.

Figure~\ref{figwavelet} shows wavelet--vaguelette deconvolution
estimators~(\ref{hfcmfixed}) obtained for $\alpha= 0$, $\alpha= 1$ and
$\alpha= 2$. Note that for moderate values of $\alpha$, the
wavelet--vaguelette estimator adjusts to spatially inhomogeneous noise
quite well. Indeed, fluctuations at the right end of the graph appear
even when $\alpha=0$ (upper left) and are due to the relatively crude
choice of threshold in formula (\ref{hfcmfixed}). Actually, for
$\alpha=0$ the noise cease to be inhomogeneous and estimator
(\ref{hfcmfixed}) reduces to Fourier-wavelet estimator of
\citet{jpkr}.
%
%
\begin{figure}

\includegraphics{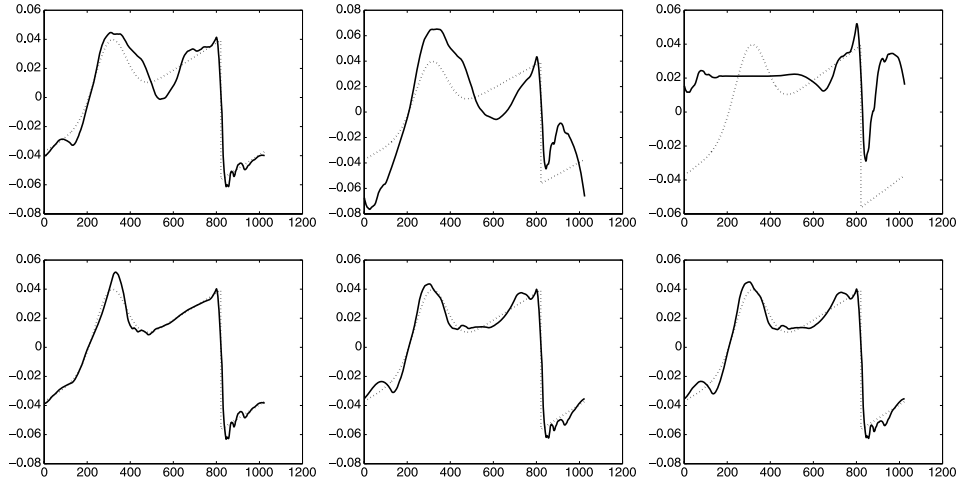}

\caption{Wavelet--vaguelette (top row) and hybrid (bottom row)
estimators for various values of $\alpha$:
$\alpha=2.5$ (left), $\alpha=3$ (middle) and $\alpha=4$ (right).
True $f$ (dotted line), wavelet--vaguelette estimator or hybrid
estimator (solid line).}\label{figvaralpha}
\end{figure}

As the values of $\alpha$ grow, SNR declines and the wavelet--vaguelette
estimators (\ref{hfcmfixed}) deteriorate. If $\alpha=3$, as
Figure~\ref{figvaralpha} shows, the wavelet--vaguelette reconstruction has
little resemblance to the regression function which it estimates. For
large values of $\alpha$, we construct hybrid estimators described in
Sections~\ref{secestimationstrategies} and \ref{secadaptive}.
Construction of the adaptive hybrid estimator consists of the following
steps.
\begin{enumerate}[3.]
\item Fix the lowest resolution level $m_1$ and the highest
resolution level $J$. For each value of $m = m_1,\ldots, J-1$,
repeat steps 2--6.

\item Obtain the wavelet--vaguelette estimator of $f$ with the
lowest resolution level $m$ using formula (\ref{hfcmfixed}).

\item Identify sets $K_{0m}$, $K_{0m}^c$, $K_{1j}$ and $K_{1j}^c$
for $j=m,\ldots, J-1$. Also, find set $\Omega_m$.

\item Form matrices $\mathbf{A}^{(m)}$ and $\mathbf{B}^{(m)}$ and
vector $\hat{\mathbf{c}}{}^{(m)}$ using formulae (\ref{matrAB})
and (\ref{vecchc}), respectively, and obtain solution
$\hat{\mathbf{z}}{}^{(m)}$ of the system of equations
(\ref{eqapproxsys}). Use vector $\hat{\mathbf{z}}{}^{(m)}$ as
coefficients $\hat{a}_{mk}$, $k \in K_{0m}$, in the
zero-affected portion of the estimator (\ref{hfom}).

\item Replace estimators of the scaling coefficients (if $k \in
K_{0m}$) and wavelet coefficients (if $k \in K_{1j}$, $j =
m,\ldots, J-1$) by zeros to obtain the zero-free portion of the
estimator (\ref{hfcmfixed}).

\item Combine wavelet coefficients in steps 4 and 5 to obtain
wavelet coefficients of~$\hat{f}_{m}$. Recover $\hat{f}_{m}$
using inverse wavelet transform. Set
\[
\lambda_m^{-2} = \sigma^2 \sum
_k \mu^{-2} \bigl(2^{-m} k \bigr) \|
T_{mk}\|^2.
\]

\item For each $m = m_1,\ldots, J-1$ and $j = m+1,\ldots, J-1$,
    evaluate matrix of the adjusted differences
%
%
\begin{equation}
\label{eqLepmatr}
\mathcal{L}_{mj} = \bigl\| ( \hat{f}_{m}-
\hat{f}_{j}) \mathbb{I}(\Omega_m) \bigr\|^2 / \bigl(
\sigma^2 n^{-1} \log n \lambda_j^{-2}
\bigr).
\end{equation}
\end{enumerate}
Choose $\hat{m}$ in (\ref{mopt}) by comparing maximum value of row $m$
of matrix $\mathcal{L}$ with a~constant $\kappa^2$
%
%
\begin{equation}
\label{moptex} \hat{m}= \min \bigl\{m\dvtx\mathcal{L}_{mj} \leq
\kappa^2 \mbox{ for all } j, m \leq j \leq J-1 \bigr\}.
\end{equation}
%

In our simulations, we used $n=1024$, $m_1=1$ and $J=7$ and carried out
hybrid estimation with $\alpha=2.5$, $\alpha=3$ and $\alpha=4$.
Simulation results for these three cases are presented in Figure~\ref{figvaralpha}. The upper and the lower rows present reconstructions
of $f$ by wavelet--vaguelette and hybrid estimator, respectively.
Observe that for \mbox{$\alpha= 2.5$} the wavelet--vaguelette estimator still
generally follows the true function $f$ but for $\alpha= 3$ or $\alpha=
4$ it bears little resemblance to $f$. The hybrid estimator allows to
account for inhomogeneity of the noise and to significantly improve
reconstruction of $f$.
%
%
\begin{figure}

\includegraphics{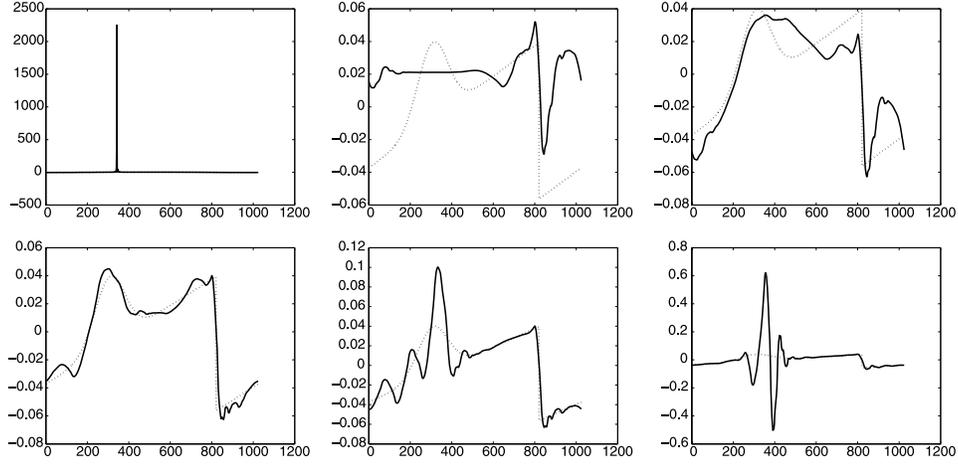}

\caption{Hybrid estimation with $\alpha=4$. True $H$ (dotted line) and
observed data (solid line) (upper left),
true $f$ (dotted line) and wavelet--vaguelette estimator (solid line)
(upper middle), true $f$ (dotted line) and
hybrid estimators (solid line) with $m=2$ (upper right), $m=3$ (lower
left), $m=4$ (lower middle)
and $m=5$ (lower right). Lepski method selects estimator with $\hat{m}
=3$ (lower left).}\label{figalpha400}
\end{figure}
%
%
\begin{figure}[b]

\includegraphics{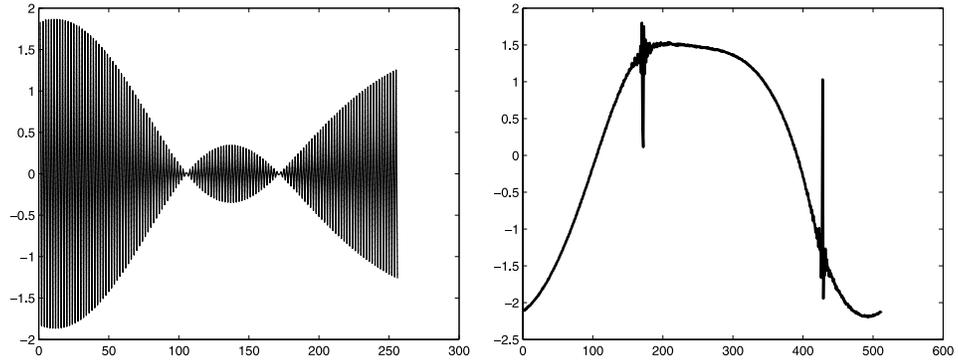}

\caption{Observed values of the signal and the true signal.
Left: the left half (first $n/2=256$ points) of the true signal (solid
line) and observed data with homogeneous noise (dashed line).
Right: true signal divided by $\mu$ (dotted line) and observed data
with heteroscedastic noise (solid line).
Here, noise level $\sigma= 0.01$, $\theta= \pi/6$, $x_{01} = 1/3$ and
$x_{02}= 5/6$.}\label{figAMsignal}
\end{figure}

Lepski procedure provides a choice of resolution level $\hat{m}$ for
each of the values of $\alpha$. Figure~\ref{figalpha400} demonstrates
hybrid estimators with the various lowest resolution levels $m$ when
$\alpha=4$. In this case, maximum over $j$ of $\mathcal{L}_{mj}$ is
very large for $m\leq2$ and is below 3 for $m \geq3$, so that
$\hat{m}=3$. Figure~\ref{figalpha400} confirms that Lepski procedure
makes the correct choice.

\subsection{Real data application}\label{secrealdata}

Below, we consider application of the hybrid estimator developed in the
paper to recovery of a convolution signal transmitted via amplitude
modulation described in Example \ref{ex3}. Mathematically, the problem
reduces to deconvolution with a spatially inhomogeneous kernel in
Section~\ref{secdeconv} and appears in the form of equation
(\ref{eqheteroscedastic}) with $\mu(x) = \cos(2\pi\omega x + \theta)$,
where $\omega\approx n/2$ and $\theta\in[0; \pi]$.
We chose $\omega= n/2 +1/2$, so that $\mu(x)$ has two zeros of order
$\alpha=2$, $x_{01} \leq1/2$ and $x_{02} = x_{01} + 1/2$ in $[0,1]$. In
particular, $x_{01} = 0.5 - \theta/\pi$ if $\theta< \pi/2$ and $x_{01}
= 1.0 - \theta/\pi$ if $\pi/2 \leq\theta\leq\pi$.

For simplicity, we considered the same set up as in simulation example,
that is, we used $q(x)=q_1 (x)$ with $\lambda= 5$ where $q_1 (x)$ is
defined in (\ref{eqvariousq}) and one of the standard test functions,
\texttt{blip}, as the true function $f(x)$. We carried out simulations
with $n = 512$, $\sigma= 0.01$ and degree 8 Daubechies wavelets. The
locations of zeros were estimated from the data. For example, when
$\theta= \pi/6$, $x_{01} = 1/3$ and $x_{02}= 5/6$, locations of zeros
were estimated as $\hat{x}_{01} = 0.33496$ and $\hat{x}_{02} =
0.83496$.

Figure~\ref{figAMsignal} presents signal $y$ with uniform noise
generated according to equation~(\ref{eqheteroscedastic}) as well as
signals $Y$ with heteroscedastic noise obtained according to equation
(\ref{eqdifnoise}) by dividing equation (\ref{eqheteroscedastic}) by
$\mu(i/n)$. Due to the limited resolution of the pictures, in
Figure~\ref{figAMsignal}, we plotted only the first $n/2 = 256$ values
of the original signals $y$.

\begin{figure}

\includegraphics{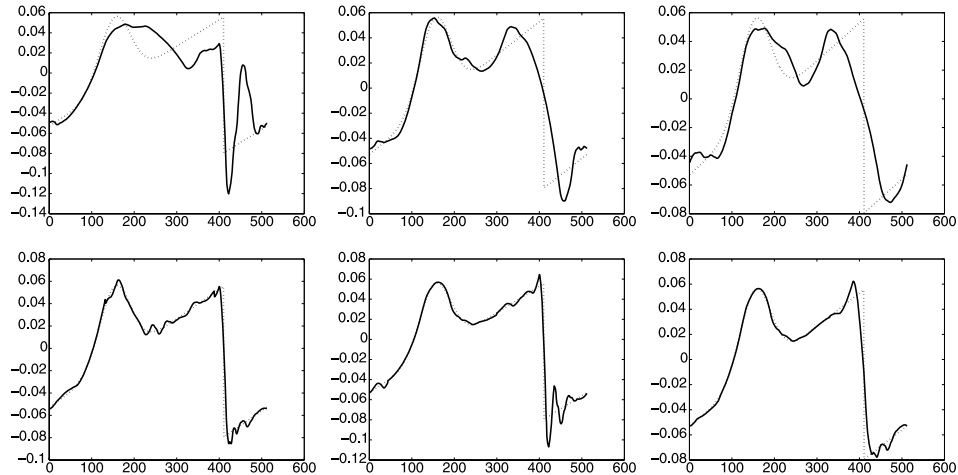}

\caption{Estimators of the signal. Top: true function $f$ (dotted
line) and wavelet--vaguelette estimator (solid line).
Bottom: true function $f$ (dotted line) and hybrid estimator (solid line).
Here, $\theta= 0.49\pi$ ($x_{01} = 0.010$ and $x_{02} = 0.510$) (left),
$\theta= \pi/6$ ($x_{01} = 1/3$ and $x_{02} = 5/6$) (middle) and
$\theta= \pi/10$ ($x_{01} = 0.400$ and $x_{02} = 0.900$) (right);
noise level $\sigma= 0.01$.}\label{figAMestimators}
\end{figure}

Figure~\ref{figAMestimators} shows wavelet--vaguelette deconvolution
estimators and hybrid estimators for various values of $\theta$:
$\theta= 0.49 \pi$ ($x_{01} = 0.010$ and $x_{02} = 0.510$), $\theta=
\pi/6$ ($x_{01} = 1/3$ and $x_{02} = 5/6$) and $\theta= \pi/10$
($x_{01} = 0.400$ and $x_{02} = 0.900$). The upper and the lower rows
present the wavelet--vaguelette and the hybrid estimators, respectively.
In all three cases, the wavelet--vaguelette estimators deteriorate due
to high value of $\alpha$ while hybrid estimators deliver satisfactory
recovery of the true function $f$. Note that the hybrid estimator has
the worst performance when $\theta= \pi/6$ due to proximity of zeros,
$x_{01} = 1/3$ and $x_{02} = 5/6$, of $\mu(x)$ to discontinuities in
the derivatives of the underlying signal $f$.

%
%


\section{Discussion}
\label{secdiscussion}

In the present paper, we consider estimation of a solution of a
spatially inhomogeneous linear inverse problem (\ref{eqmaineq}) with
possible singularities. The special feature of problems like this is
that the degree of ill-posedness depends not only on the scale but also
on location. In spite of a huge number of publications devoted to
linear inverse problems, to the best of our knowledge, this type of
problems has never been treated before.
We consider a version of a spatially inhomogeneous problem where there
exists a singularity point $x_0$ such that the norm of the solution
grows when the right-hand side is localized in the vicinity of $x_0$.
%
We characterize ill-posedness and spatial inhomogeneity of operator $Q$
in terms of wavelet--vaguelette decomposition. The novel feature here is
that the norms of vaguelettes depend on location and may be infinite in
the vicinity of a singularity point, so that SVD-type solutions cease
to work.

For this reason, estimators obtained in the paper are based either on
wavelet--vaguelette decomposition (if the norms of all vaguelettes are
finite) or on a hybrid of wavelet--vaguelette decomposition and Galerkin
method (if vaguelettes in the neighborhood of the singularity point
have infinite norms).
We show that, up to a logarithmic factor, the hybrid estimator
attains the asymptotically optimal 
convergence rates.

The theory presented in the paper is supplemented by examples of
deconvolution with a spatially inhomogeneous kernel, deconvolution in
the presence of locally extreme noise or extremely inhomogeneous
design.
In addition, we apply the technique to recovery of a convolution signal
transmitted via amplitude modulation.

We note that the wavelet-based estimation procedure presented in the
paper is motivated by the need of constructing an asymptotically
optimal estimator in the case when the unknown function $f$ is
spatially inhomogeneous. The estimator uses relatively crude
thresholding procedure which can be improved by applying more sensitive
thresholding techniques. Moreover, one can possibly find more efficient
computational procedures than the hybrid estimator if establishing
asymptotic optimality is not a priority.

The paper assumes that the operator $Q$ in (\ref{eqmaineq}) is
completely known. However, if this is not true, it will be interesting
to investigate how uncertainty about $Q$ affects the rates of
convergence. Also, although the hybrid estimator works adequately when
$Q$ is completely known, it would require appropriate modifications if
$Q$ is partially unknown.

Finally, in the paper, we consider only the simplest case when the
unknown function is univariate and is defined on an interval. The
problem can be naturally extended to the case of multivariate function
$f$ which belongs to an isotropic or anisotropic Besov space. However,
all these extensions will be a matter of future investigation.


\section*{Acknowledgements}

The author would like to thank the Editor, the Associate Editor and two
anonymous referees for their suggestions and support that helped to
significantly improve the paper



\begin{supplement}[id=suppA]
\stitle{Proofs}
\slink[doi]{10.1214/13-AOS1166SUPP}
\sdatatype{.pdf}
\sfilename{aos1166\_supp.pdf}
\sdescription{Supplement contains proofs of the statements in the manuscript.}
\end{supplement}

\printaddresses

\end{document}